\documentclass[11pt,twoside]{amsart}

\usepackage[colorlinks=true, linkcolor=blue, anchorcolor=blue, citecolor=blue, filecolor=blue, menucolor= blue, pagecolor=black, urlcolor=red]{hyperref}

\newtheorem{thm}{Theorem}[section]
\newtheorem{prop}[thm]{Proposition}
\newtheorem{lem}[thm]{Lemma}

\newtheorem{defn}[thm]{Definition}

\newtheorem{example}[thm]{Example}
\newtheorem{rem}[thm]{Remark}

\newcommand{\pr}[1]{\mathbb{P}^{#1}}
\newcommand{\skipit}[1]{{}}
\newcommand{\prfend}{\hbox to7pt{\hfil}
\par\vskip-\baselineskip\hbox to\hsize
{\hfil\vbox {\hrule width6pt height6pt}}\vskip\baselineskip}

\newcommand{\Oc}{\mathcal{O}}
\DeclareMathOperator{\cok}{cok}
\DeclareMathOperator{\Neg}{Neg}
\DeclareMathOperator{\nneg}{neg}
\DeclareMathOperator{\Cl}{Cl}
\DeclareMathOperator{\Image}{Im}
\DeclareMathOperator{\EFF}{EFF}
\DeclareMathOperator{\NEF}{NEF}
\newcommand {\C}[1]{\mathcal{#1}}
\DeclareMathOperator{\Pic}{Pic}
\DeclareMathOperator{\Char}{char}

\newcommand{\myarrow}[2]{\hbox to #1pt{\hfil$\to$\hfil}{\hskip-#1pt{\raise
10pt\hbox to#1pt{\hfil$\scriptscriptstyle #2$\hfil}}}}

\begin{document}

\title{Classifying Hilbert functions of fat point subschemes in $\pr{2}$}

\author{A.\ V.\ Geramita, B.\ Harbourne \& J.\ Migliore}

\address{A.\ V.\ Geramita\\
Department of Mathematics\\
Queen's University\\
Kingston, Ontario, and\\
Dipartimento di Matematica\\ 
Universit\`a di Genova\\ 
Genova, Italia}
\email{anthony.geramita@gmail.com}

\address{Brian Harbourne\\
Department of Mathematics\\
University of Nebraska\\
Lincoln, NE}
\email{bharbour@math.unl.edu}

\address{Juan Migliore\\
Department of Mathematics\\
University of Notre Dame\\
South Bend, IN}
\email{migliore.1@nd.edu}

\date{November 23, 2010}

\markboth{A.\ V.\ Geramita, B.\ Harbourne, J.\ Migliore}{Classifying 
Hilbert functions of fat point subschemes in $\pr{2}$}

\thanks{Acknowledgments: We thank the NSF, whose 
support for the MAGIC05 conference at the University of Notre Dame gave
us an opportunity to begin work on this paper. We thank Zach Teitler for 
pointing out a typo, and we especially thank the referee for being so
prompt and so thorough.
We also thank Margaret Bayer, Jonathan Cutler, Jeremy 
Martin, Jamie Radcliffe, Hal Schenck, and Peter Vamos 
for their comments and suggestions regarding matroids. Geramita thanks 
the NSERC for research support, while Harbourne and Migliore thank the NSA for
supporting their research
(most recently, under grants H98230-07-1-0066
and H98230-07-1-0036, respectively).}

\begin{abstract}
The paper \cite{GMS} raised the question of what the possible
Hilbert functions are for fat point subschemes 
of the form $2p_1+\cdots+2p_r$, for all possible choices
of $r$ distinct points in $\pr2$.
We study this problem for $r$ points in $\pr2$
over an algebraically closed field $k$ of arbitrary characteristic
in case either $r\le8$ or the points lie on a (possibly reducible) conic. 
In either case, it follows from \cite{freeres} and \cite{BHProc}
that there are only finitely many configuration types of points,
where our notion of configuration type is a generalization of the notion of a
representable combinatorial geometry, also known as a 
representable simple matroid. (We say
$p_1,\ldots,p_r$ and $p'_1,\ldots,p'_r$ have the same {\it configuration type\/}
if for all choices of nonnegative integers $m_i$,
$Z=m_1p_1+\cdots+m_rp_r$ and $Z'=m_1p'_1+\cdots+m_rp'_r$
have the same Hilbert function.)  
Assuming either that $7\le r\le 8$ (see \cite{GuH} for the cases $r\le 6$) or
that the points $p_i$ lie on a conic, we explicitly
determine all the configuration types, and show how the
configuration type and the coefficients $m_i$
determine (in an explicitly computable way)
the Hilbert function  (and sometimes the graded Betti numbers) of
$Z=m_1p_1+\cdots+m_rp_r$. 
We demonstrate our results by explicitly listing all
Hilbert functions for schemes of $r\le 8$ double points, and for each Hilbert function
we state precisely how the points must be arranged (in terms of
the configuration type) to obtain that Hilbert function.
\end{abstract}



\maketitle

\section{Introduction}\label{intro}

Macaulay completely solved the problem of describing what can be the Hilbert function of a homogeneous ideal $I$ in a polynomial ring $R$.  His solution, however, left open the question of which of these functions are the Hilbert functions of ideals in specific classes, such as prime homogeneous ideals (i.e., ideals of irreducible projective varieties) or, what will be a  focus in this paper, symbolic powers of ideals of reduced zero-dimensional subschemes of projective 2-space.   

Various versions of these questions have been studied before.
A complete determination of the Hilbert functions of reduced
0-dimensional subschemes of projective space is given in \cite{GMR}.  But
it is not, for example, known which functions arise as Hilbert functions  
of points taken with higher multiplicity, even if the induced reduced subscheme 
consists of generic points; see, for example, 
\cite{Ci} and its bibliography, \cite{Gi}, \cite{Vanc}, \cite{Hi} and \cite{Ro}.
The paper \cite{GMS}
asks what can be said about Hilbert functions and graded Betti numbers
of symbolic powers $I^{(m)}$ of the ideal $I$ 
of a finite set of reduced points in $\pr2$, particularly when $m = 2$,
while the papers \cite{GuV1} and \cite{GuV2} study $I^{(m)}$ for larger $m$,
but only when $I$ defines a reduced subscheme of projective 
space whose support is close to a complete intersection. 
Other work has focused on what one knows as a consequence of knowing 
the Hilbert function; e.g., \cite{BGM} shows how the growth of 
the Hilbert function of a set of points influences its geometry,
\cite{Cam} studies how the Hilbert function
constrains the graded Betti numbers in case $I$ has height 2, and
\cite{ER} studies graded Betti numbers but more generally 
for graded modules $M$ over $R$. 

In this paper we will determine precisely what functions occur as
Hilbert functions of ideals $I(Z)$ of fat point subschemes $Z$
of $\pr2$ over an algebraically closed field $k$ of arbitrary characteristic, 
in case either $I(Z)$ contains some power of a form of degree 2 
(i.e., the support of $Z$ lies on a conic), or in case the support of $Z$ consists
of 8 or fewer points. If either the support of $Z$ lies on a conic or consists
of 6 or fewer points, by results of \cite{Cat}, \cite{freeres}, and \cite{GuH},
one can also determine the graded Betti numbers of $I(Z)$.
In the case of ideals 
of the form $I(Z)$ when the support of $Z$ is small or lies on
a conic, our results allow us to give explicit answers for many of the questions raised in the papers cited above.
For example, given the Hilbert function $h$ of $I(Z')$ for a reduced subscheme $Z'=p_1+\cdots+p_r$,
we can explicitly determine the Hilbert functions of all symbolic powers $I(Z)^{(m)}$
for all $Z$ whose ideal $I=I(Z)$ has
Hilbert function $h$, as long as $Z'$ has support on a conic or at 8 or fewer points,
and we can in addition explicitly determine the graded Betti numbers
of $I(Z)^{(m)}$ 
as long as $Z'$ has support either on a conic (by \cite{Cat}, \cite{freeres}) or at 6 or fewer points
(by \cite{GuH}).

We now discuss the context of our work in more detail.
Let $Z$ be a subscheme of $\pr2$. Let 
$I(Z)$ be the corresponding saturated homogeneous ideal in the polynomial ring
$R=k[x,y,z]$ (over an algebraically closed field $k$ of arbitrary characteristic) 
which defines $Z$. Then $R(Z)=R/I(Z) = \oplus_{j \geq 0} R(Z)_j$ is a graded
ring whose Hilbert function $h_Z$  is defined by $h_Z(t)=\hbox{dim}_kR(Z)_t$.

If $Z$ is a fat point subscheme (defined below), then $h_Z$ is a nondecreasing
function with $h_Z(0)=1$ and such that $h_Z(t)=\hbox{deg}(Z)$
for all $t\ge\hbox{deg}(Z)-1$. Thus there are only finitely many 
possible Hilbert functions for a fat point subscheme $Z$ of given degree.
In case $Z$ is smooth (i.e. a finite set of distinct points) it is known what these finitely many
Hilbert functions are \cite{GMR}. The paper \cite{GMS} raises the problem of 
determining the Hilbert functions, and even the graded Betti numbers
for a minimal free resolution of $R(Z)$ over $R$, for fat point
subschemes of $\pr2$ more generally.

If $Z$ is a reduced 0-dimensional subscheme of $\pr2$, let $mZ$
denote the subscheme defined by the $m$th symbolic
power $I(mZ)=(I(Z))^{(m)}$ of the ideal $I(Z)$.
If $Z$ consists of the points $p_1,\ldots, p_r$, we write
$Z=p_1+\cdots+p_r$, and we have $I(Z)=I(p_1)\cap \cdots\cap I(p_r)$,
where $I(p_i)$ is the ideal generated by all forms of $R$ that vanish
at $p_i$. We also write $mZ=mp_1+\cdots+mp_r$.
The ideal $I(mZ)$ is just $I(p_1)^m\cap \cdots\cap I(p_r)^m$.
More generally, given any nonnegative integers $m_i$,
the {\it fat point\/} subscheme denoted $m_1p_1+\cdots+m_rp_r$
is the subscheme defined by the ideal 
$I(m_1p_1+\cdots+m_rp_r)=I(p_1)^{m_1}\cap \cdots\cap I(p_r)^{m_r}$.
Given a fat point subscheme $Z=m_1p_1+\cdots+m_rp_r$, we refer to the 
sum of the points for which $m_i>0$ as the {\it support\/} of $Z$ and we define the {\it degree} of $Z$, denoted $\deg Z$, to be the integer $\sum_{i=1}^r \binom{m_i +1}{2}$.

Given the Hilbert function $h$ of a reduced fat point subscheme, 
the focus of the paper \cite{GMS} is to determine what
Hilbert functions occur for $2Z$, among all reduced $Z$ whose 
Hilbert function is $h$. The paper \cite{GMS} represents
a step in the direction of answering this problem:
For each $h$ which occurs for a reduced fat point subscheme of $\pr2$,
\cite{GMS} determines the Hilbert function of $2Z$
for an explicitly constructed subscheme $2Z$ whose support
$Z$ has Hilbert function $h$. This leaves open the problem of
giving a complete determination of the Hilbert functions which occur
among all fat points subschemes $2Z$ whose support $Z$
has Hilbert function $h$, but \cite{GMS} proves
that there is a maximal such Hilbert function.
This proof is nonconstructive; what this
maximal Hilbert function is and how it can be found 
is an open problem. The paper \cite{GMS} also raises the question
of whether there is a minimal Hilbert function
(among Hilbert functions for all $2Z$ such that
$Z$ is reduced with given Hilbert function). 

In this paper, we give a complete answer to all
of the questions raised in \cite{GMS} about Hilbert functions
if the degree of the support is
8 or less or the support is contained in a conic.  
We classify all possible arrangements of $r$ points $p_i$
in case either $r\le8$ or the points are contained in a conic,
and determine the Hilbert functions for all fat point subschemes
$Z=m_1p_1+\cdots+m_rp_r$, regardless of the multiplicities.
This is possible because Bezout considerations (see \cite{freeres},
\cite{BHProc} and Example \ref{ptsonconicex})
completely determine the Hilbert functions for any $Z$ if either $r\le 8$ or
the points $p_i$ lie on a conic (and these same considerations
even determine the graded Betti numbers of $I(Z)$
if either $r\le 6$ or the points lie on a conic).

Without a bound on the degree of $Z$,
there are clearly infinitely many Hilbert functions. 
However, results of \cite{freeres} (in case the points
$p_1,\dots,p_r$ lie on a conic) or of \cite{BHProc} (in case 
$r\le 8$) imply that to each set of points $p_1,\dots,p_r$
we can attach one of a finite set of combinatorial structures
we call {\it formal configuration types\/}. Those
formal configuration types which are representable
(meaning they arise from an actual set of points)
turn out to be the same thing as what we define to be a configuration type.
In case $r\le8$ or the points $p_1,\ldots,p_r$ lie on a conic,
configuration types have the property that if we know the integers $m_i$
and the configuration type of the points $p_i$, we can 
explicitly write down the Hilbert function
of $m_1p_1+\cdots+m_rp_r$ (and, if either $r\le 6$ or the
points lie on a conic, we can even write down the graded Betti numbers).

The ultimate goal of the program started in \cite{GMS} is to determine
all Hilbert functions (and the graded Betti numbers) for
fat point subschemes $Z$ whose support has given
Hilbert function $h$. Note that $h(2)$ is less than 6 if and only if the points of the support
lie on a conic. Thus in case $h(2)<6$, the results of \cite{freeres}
give a complete answer to the program of \cite{GMS}, once one writes
down the configuration types of points on a conic (which is easy).
The results of \cite{BHProc} similarly
give a complete answer to the Hilbert function part of
the program of \cite{GMS} when $h$ is less than or equal to the constant function 8,
modulo writing down the configuration types.
But unlike the case of points on a conic, determining the configuration types
for 7 or 8 distinct points takes a fair amount of effort. That determination is given 
here for the first time.

Giving the list of types for 8 points subsumes doing so for any fewer number of points,
but it is of interest to consider the case of 7 points explicitly.
The case of 6 points is worked out in \cite{GuH}. (In fact, the main
point of \cite{GuH} is to determine graded Betti numbers. Thus
\cite{GuH} entirely resolves the program of \cite{GMS}
in case $h$ is less than or equal to the constant function 6. 
The case of fewer than 6 points
can easily be recovered from the case of 6 points.)
It turns out that there are 11 different configuration types
of 6 distinct points of $\pr2$.
We find 29 types for 7 points and 143 types for 8 points.

However, a new feature arises for the case of 7 or 8 points,
compared with 6 points.
Formal configuration types, as we define them, are matroid-like combinatorial objects
which can be written down without regard to whether or not
some set of points exists having that type.
A combinatorial geometry (or, equivalently, a simple matroid) 
whose points have span of dimension $s\le 2$
is a matroid, of rank $s+1\le 3$, without loops or parallel
elements. Our classifications of formal configuration types of
sets of 7 (respectively 8) distinct points
includes the classification of all combinatorial geometries on 
7 (respectively 8) points whose span has dimension at most 2 \cite{BCH}
(alternatively, see \cite{MR}).
A key question for combinatorial geometries is that of representability;
that is: does there exist a field for which there occurs an actual set of
points having the given combinatorial geometry?
This is also an issue for our configuration types.
In the case of 6 points, every formal configuration type is representable for
every algebraically closed field, regardless of the characteristic.
For 7 points, every formal configuration type is representable for some
algebraically closed field, but sometimes representability
depends on the characteristic. 
For 8 points, there are three formal configuration types 
which are not representable at all.

Given only the configuration type and the multiplicities $m_i$,
we also describe an explicit procedure (which, if so desired,
can be carried out by hand) for computing the Hilbert function of 
any $Z=m_1p_1+\cdots+m_rp_r$ when either $r\le 8$ or
the points lie on a conic. If the points lie on a conic (using \cite{Cat} or \cite{freeres}) or 
if $r\le 6$ (using
\cite{GuH}), one can also determine the graded Betti numbers of the
minimal free resolution of $R(Z)$.
Example \ref{ptsonconicex} demonstrates the procedure
in detail in case the points lie on a conic.
In case $6\le r\le 8$,
the procedure is implemented as a web form that can be run
from any browser. To do so, visit
\url{http://www.math.unl.edu/~bharbourne1/FatPointAlgorithms.html}.

\section{Background}\label{bckgrnd}

In this section we give some definitions
and recall well known facts that we will need later.
First, the minimal free resolution of $R(Z)=R/I(Z)$ for a fat point subscheme $Z\subset \pr2$
is
$$0\to F_1\to F_0\to R\to R(Z)\to 0$$
where $F_0$ and $F_1$ are free graded $R$-modules
of the form $F_0=\oplus_{i\ge 0} R[-i]^{t_i}$ and
$F_1=\oplus_{i\ge 0} R[-i]^{s_i}$. The indexed values
$t_i$ and $s_i$ are the graded Betti numbers of $Z$.
Consider the map 
$$\mu_i:I(Z)_i\otimes R_1\to I(Z)_{i+1}$$
given by multiplication of forms of degree $i$ in $I(Z)$
by linear forms in $R$. Then $t_{i+1}=\dim\cok(\mu_i)$,
and $s_i=t_i+\Delta^3h_Z(i)$ for $i\ge 1$, where $\Delta$ is the difference operator
(so $\Delta h_Z(i+1)=h_Z(i+1)-h_Z(i)$).
Note that knowing $t_i$ and the Hilbert function of $Z$
suffice to ensure that we know $s_i$.
Thus once $h_Z$ is known, if we want to know the graded Betti numbers
for $Z$, it is enough to determine $t_i$ for all $i$.

We now introduce the notion of Hilbert function equivalence for  
ordered sets of points in $\pr2$. 

\begin{defn}\rm Let $p_1,\ldots,p_r$ and $p'_1,\ldots,p'_r$ be ordered 
sets of distinct points of $\pr2$. We say these sets are 
{\it Hilbert function equivalent\/} if
$h_Z=h_{Z'}$ for every choice of nonnegative integers $m_i$,
where $Z=m_1p_1+\cdots+m_rp_r$ and $Z'=m_1p'_1+\cdots+m_rp'_r$.
We refer to a Hilbert function equivalence class as a 
{\it configuration type}.
\end{defn}

\begin{rem}\label{howmany}\rm Given $r$ and $m_1, \ldots , m_r$ it is clear that there are only finitely many possible Hilbert functions for fat point schemes $Z = m_1p_1+\cdots+m_rp_r$.  However, for any given $r$ there is absolutely no guarantee that there are finitely many configuration types for sets of $r$ points. Indeed, there are infinitely many types
for $r=9$ points (see Example \ref{infmanytypes}). Thus it is interesting that 
by our results below there are only finitely many configuration types among
ordered sets of $r$ distinct points if either the points lie on a conic or $r \leq 8$ points.  This allows us to answer many of the questions raised in \cite{GMS} (and more) for these special sets of points.  
\end{rem}

The methods used in this paper are very different from those used in \cite{GMS}, partly because the focus in this paper is on special sets of points, in particular points on a conic and sets of $r \leq 8$ points in $\pr2$. For these cases, the second author has demonstrated, in a series of papers, the efficacy of using the theory of rational surfaces, in particular the blow-ups of $\pr2$ at such points.

The tools of this technique are not so familiar in the commutative algebra community (where the original questions were raised) and so we thought to take this opportunity not only to resolve the problems of \cite{GMS} for such sets of points, but also to lay out the special features of the theory of rational surfaces which have proved to be so useful in these Hilbert function investigations.

We begin with basic terminology and notation.  Let $p_1,\ldots,p_r$ be distinct points of $\pr2$.
Let $\pi:X\to \pr2$ be the morphism obtained by blowing up the points.
Then $X$ is a smooth projective rational surface.
Its divisor class group $\Cl(X)$ is a free Abelian group
with basis $L$, $E_1, \ldots,E_r$, where $L$ is the class
of the pullback via $\pi^*$ of the class of a line, and
$E_i$ is the class of the fiber $\pi^{-1}(p_i)$.
We call such a basis an {\it exceptional configuration\/} for $X$.
It is an orthogonal basis for $\Cl(X)$ with respect to the
bilinear form given by intersecting divisor classes.
In particular, $-L^2=E_i^2=-1$ for all $i$, and
$L\cdot E_i=0=E_i\cdot E_j$ for all $i$ and all $j\ne i$.
An important class is the {\it canonical class} $K_X$.
In terms of the basis above it is $K_X=-3L+E_1+\cdots+E_r$.
Given a divisor $F$, we write
$h^0(X, \Oc_X(F))$ for the dimension of the space
$H^0(X, \Oc_X(F))$ of global sections of the line bundle
$\Oc_X(F)$ associated to $F$.
To simplify notation, we will hereafter
write $h^0(X, F)$ in place of $h^0(X, \Oc_X(F))$, etc.  
Since $h^0(X, F)$ is the same for all divisors $F$ in the same divisor class,
we will abuse notation and use $h^0(X, F)$ even when $F$ is just a divisor class. 
The projective space ${\mathbb P} ( H^0(X, \Oc_X(F))$ will be denoted $|F|$.

Recall that a {\it prime divisor} is the class of a reduced irreducible curve
on $X$, and an {\it effective divisor}  is a nonnegative integer
combination of prime divisors. For simplicity, we will refer to a 
divisor class as being {\it effective\ } if it is the class of an effective divisor.  We will denote the collection of those effective divisor classes by EFF$(X)$.  Notice that for the cases we are considering in this paper (i.e. $r$ points on a conic or $r \leq 8$ points) the {\it anticanonical divisor class} $-K_X$ is effective.  We say that a divisor $F$ is {\it ample} if $F\cdot D > 0$ for every effective divisor $D$.

Our interest in this paper is to decide when a divisor class $F$ on $X$ is effective and then to calculate $h^0(X,F)$.  To see why this is relevant to the problem of computing Hilbert functions
of ideals of fat points, let $Z=m_1p_1+\cdots+m_rp_r$. Then under the identification
of $R_i$ with $H^0(\pr2, \Oc_{\pr2}(i))=H^0(X, iL)$, we have $I(Z)_i=H^0(X, F(Z,i))$, where
$F(Z,i)=iL-m_1E_1-\cdots-m_rE_r$. We also have a canonical map
$$\mu_F: H^0(X, F)\otimes H^0(X, L)\to H^0(X, F+L)$$
for any class $F$. In case $F=F(Z,i)$, then $t_{i+1}=\dim\cok(\mu_{F(Z,i)})$.  As a result, we can use facts about divisors on $X$, in particular about $h^0(X,F)$, 
to obtain results about the Hilbert function and graded Betti numbers for $Z$.

Since $F$ is a divisor on a rational surface, the Riemann-Roch Theorem gives
$$
h^0(X,F) - h^1(X,F) + h^2(X,F) = \frac{F\cdot F - K_\cdot F}{2} + 1
$$
So, if for some reason we were able to prove that the divisor $F$ we are interested in has $h^1(X,F)=h^2(X,F) = 0$ then the Riemann-Roch formula would calculate $h^0(X,F)$ for us.

A divisor whose intersection product with {\bf every} effective divisor is $\geq 0$ is called {\it numerically effective (nef)}. The collection of nef divisor classes forms a cone which is denoted NEF$(X)$.  It is always true for a nef divisor $F$ that $h^2(X, F)=0$ (since $(K-F) \cdot L<0$, we see that $K-F$ is not effective; now use Serre duality).  So, if we could somehow connect the divisor we are interested in to a nef divisor, then at least we might be able to get rid of the need to calculate $h^2(X,F)$.  Unfortunately, it is not true, in general, that $h^1(X, F)=0$ for $F$ a nef divisor. 

However, one of the fundamental results in this area asserts that for rational surfaces obtained by blowing up either any number of points on a conic, or any $r \leq 8$ points, nef divisors do have $h^1(X,F) = 0$.  More precisely, 

\begin{thm}\label{BettisEtc} 
Let $X$ be obtained by blowing up $r$ distinct points of $\pr2$.
Suppose $F$ is a nef divisor on $X$.  
\begin{itemize}
\item[(a)]If either $2L-E_1-\cdots-E_r$ is effective (i.e., the points $p_i$ lie on a conic), or $r\le 8$, then  $h^1(X,F) = 0 = h^2(X,F)$ and so $h^0(X,F)=\frac{(F^2-K_X\cdot F)}{2}+1$.
\item[(b)]If $2L-E_1-\cdots-E_r$ is effective, then $\mu_F$ is surjective.
\item[(c)]If $r\le6$, then $\mu_F$ has maximal rank (i.e., is either surjective or injective).
\end{itemize}
\end{thm}

\begin{proof} As indicated above, $h^2(X,F)=0$ for a nef divisor.
When $2L-E_1-\cdots-E_r$ is the class of a prime divisor,
the fact that 
$$
h^0(X,F)=\frac{(F^2-K_X\cdot F)}{2}+1
$$ 
was proved in
\cite{Vanc}, and surjectivity for $\mu_F$ was proved in \cite{Cat}.
Both results in the general case (assuming only that the points $p_i$ lie
on a conic, not that the conic is necessarily reduced and irreducible,
hence subsuming the case that the points are collinear)
were proved in \cite{freeres}.
The fact that $h^0(X,F)=(F^2-K_X\cdot F)/2+1$ if $r\le 8$ 
was proved in \cite{BHProc}. The fact that $\mu_F$
has maximal rank for $r\le 6$ is proved in \cite{GuH}.
(We note that if $r\ge 7$, it need not be true that $\mu_F$ always has maximal rank
if $F$ is nef \cite{7pts}, \cite{FHH}, and if $r=9$, it need not be true that
$h^1(X,F) = 0$ if $F$ is nef, as can be seen by taking $F=-K_X$,
when $X$ is the blow-up of the points of intersection of a pencil of cubics.) 
\end{proof} 

\begin{rem}\label{regularity}\rm   Let $Z = m_1p_1 + \cdots + m_rp_r$ and let $F = tL - m_1E_1 - \cdots - m_rE_r$ be nef.
In Theorem \ref{BettisEtc}, (a) says if either the points lie on a conic or $r\le8$ then (by means of a standard argument) $h_Z(t) = \deg Z$, i.e.  $t +1 \geq $ the Castelnuovo-Mumford regularity of $I(Z)$; $(b)$ says that $I(Z)$ has no generator in degree $t+1$ when the points lie on a conic and $(c)$ says if $r\le6$ that the number of generators of $I(Z)$ in degree $t+1$ is $\max\{ \dim I(Z)_{t+1} - 3\dim I_Z(t), 0 \}$.  Of course, $(b)$ and $(c)$ are interesting only if $t+1$ is equal to the Castelnuovo-Mumford regularity of $I(Z)$, since $\mu_t$ is always surjective for
$t$ at or past the regularity.
\end{rem} 

In general, we will want to compute $h^0(X, F)$ and the rank of $\mu_F$ even if $F$ is not nef.  The first problem is to decide whether or not $F$ is effective.  If $F$ is not effective, then $h^0(X, F)=0$ by definition, and it is clear in that case that $\mu_F$ is injective,  and hence that $\dim\cok(\mu_F)=h^0(X, F+L)$.  So, if we can decide that $F$ is not effective then we are done.  What if $F$ is effective?

Suppose $F$ is effective {\bf AND} we can find all the fixed components $N$ of $F$ (so $h^0(X, N)=1$)).  We can then write  $F=H+N$, and clearly $H$ is nef with $h^0(X,F) = h^0(X,H)$. When Theorem \ref{BettisEtc}(a) applies, this finishes our task of finding $h^0(X,F)$.  The decomposition $F = H+N$ is called the {\it Zariski decomposition} of $F$.

Moreover, the restriction exact sequence
$$
0 \rightarrow {\mathcal O}_X(-N) \rightarrow {\mathcal O}_X \rightarrow {\mathcal O}_N \rightarrow 0
$$
tensored with ${\mathcal O}_X(F + L)$, gives an injection $|H+L| \to |F+L|$ which is defined by $D \mapsto D+N$.  Now, $ \Image(\mu_H)$ defines a linear subsystem of $|H+L|$, which we'll denote $|\Image(\mu_H)|$, and the image of $\mu_F$ is precisely the inclusion of $|\Image(\mu_H)|$ in $|F+L|$.  Hence
$$
\dim\cok(\mu_F)=\dim\cok(\mu_H)+(h^0(X, F+L)-h^0(X, H+L)).
$$

So, if it were possible to find the Zariski decomposition of a divisor on $X$, then it would be possible to find $h^0(X,F)$ for any divisor class $F$ when Theorem \ref{BettisEtc} applies, since it furnishes us with most of what we need to find Hilbert functions and graded Betti numbers.

To sum up, what we need is an explicit way to decide if a divisor $F$ is effective or not, and when it is effective, a way to find its Zariski decomposition.  These are now our two main tasks.

It turns out that in order to explicitly carry out these two tasks, we have to know the {\it prime} divisors $C$ with $C^2 < 0$.  We thus define $\Neg(X)$ to be the classes of those prime divisors $C$ with $C^2<0$ and further define $\nneg(X)$ to be the subset of $\Neg(X)$ of classes of those $C$ with $C^2<-1$.

We now give the definition of the {\it negative curve type\/} of a set of points
 $p_1,\ldots,p_r\in \pr2$.

\begin{defn}\label{negcurtype} \rm Let $p_1,\ldots,p_r$ and $p'_1,\ldots,p'_r$ be
ordered sets of distinct points of $\pr2$. Let $X$ ($X'$, respectively) 
be the surfaces obtained by blowing the points up.
Let $L, E_1,\ldots,E_r$ and $L', E'_1,\ldots,E'_r$ be the respective
exceptional configurations, and let $\phi:\hbox{Cl}(X)\to\hbox{Cl}(X')$
be the homomorphism induced by mapping $L\mapsto L'$
and $E_i\mapsto E'_i$ for all $i>0$.
We say $p_1,\ldots,p_r$ and $p'_1,\ldots,p'_r$
have the same {\it negative curve type\/} if 
$\phi$ maps $\Neg(X)$ bijectively to $\Neg(X')$.
\end{defn}

What will turn out to be the case is that if either $r\le 8$ or
the points $p_1,\ldots,p_r$ lie on a conic, then 
the configuration types are precisely the negative curve types.

In Section \ref{PtsConic} we show how to carry out the two tasks mentioned above and, in the process, demonstrate that configuration type and negative curve type are the same, for the case that $X$ is obtained by blowing up $r$ points on a conic. The case of points on a line was handled in \cite{freeres} and is, in any case, subsumed by the case of points on a conic, which also subsumes the case of $r\le5$ points. In Section \ref{7and8pts} we treat the case of $7\le r\le 8$ points.

\section{Points On a Conic}\label{PtsConic}

As mentioned above, in this section we show how to carry out our two tasks for the case of any $r$  points on a conic.  This case is technically simpler than the cases of $r = 6,\ 7$, or 8 points, but has the advantage that all the major ingredients 
for handling the cases of $6, 7$ or 8 points are already present in the simpler argument for points on a conic.  What should become clear from the discussion of this special case is the reason why there are only finitely many configuration types and why they are the same as the negative curve types.

In order to state the next result efficiently, we introduce some finite families of special divisors on surfaces obtained by blowing up $\pr2$ at $r$ distinct points lying on a conic.  Let $\C B$ be the classes $\{E_1,\cdots, E_r\}$ on $X$,  obtained from the blow-ups of the points $p_1, \ldots , p_r$.  These are always classes in $\Neg(X)$ but not in $\nneg(X)$.  Let $\C L$ be the set of classes
$$
{\C L} := \{L-E_{i_1}-\cdots-E_{i_j} \mid  2\le j,\quad 0<i_1<\cdots<i_j\le r\},
$$ 
Note that these classes are all potentially in $\Neg(X)$ (such a class may not be a prime divisor; e.g., if $p_1, p_2, p_3$ are three points on a line, then $L-E_1 -E_2 = (L-E_1 -E_2 -E_3) + E_3$ is a sum of two classes, where $E_3$ is clearly prime and so is $L-E_1 -E_2 -E_3$ if there are no additional points $p_i$ on the line). If $r > 4$ then there is another effective divisor which is potentially in $\Neg(X)$, namely $Q = 2L-E_1-\cdots- E_r$.  Notice also that the divisor $A_r = (r-2)L -K_X$ meets all the divisors in the set ${\C B} \cup {\C L} \cup \{ Q \}$ positively.

There is no loss in assuming $r\ge 2$, since the case $r<2$ is easy to handle by ad hoc methods.
(For example, if $r=0$ and $F=tL$, then $F$ is effective if and only if $t\ge0$, and $F$ is nef if and only if $t\ge0$.
If $r=1$ and $F=tL-mE_1$, then $F$ is effective if and only if $t\ge0$ and $t\ge m$, and $F$ is nef if and only if
$t\ge m$ and $m\ge 0$. However, if $r\le 1$, then it is not true that every element of $\EFF(X)$ is a nonnegative integer
linear combination of elements of $\Neg(X)$.)

\begin{prop}\label{ptsonconic} 
Let $X$ be obtained by blowing up $r\ge2$ distinct points
of $\pr2$, and assume that $Q=2L-E_1-\cdots-E_r$ is effective.
\begin{itemize}
\item[(a)] Then $\NEF(X)\subset \EFF(X)$, and 
every element of $\EFF(X)$ is a nonnegative integer
linear combination of elements of $\Neg(X)$;
\item[(b)] $\Neg(X)\subset \C L \cup \C B\cup \{Q\}$;
\item[(c)] $\Neg(X)=\nneg(X) \cup \{C\in \C L \cup 
\C B\cup \{Q\} \mid C^2=-1, C\cdot D\ge 0\hbox{ for all }D\in \nneg(X)\}$;
\item[(d)] for every nef divisor $F$, $|F|$ is base point (hence fixed component)
free; and 
\item[(e)] the divisor $A_r$ is ample.
\end{itemize}
\end{prop}

\begin{proof} Parts (a, b, d) follow from Lemma 3.1.1 \cite{freeres}.
For (c), we must check that if $C\in \C L \cup \C B\cup \{Q\}$, $C^2=-1$ and 
$C\cdot D\ge 0$ for all $D\in \nneg(X)$, then $C$ is the class of a prime divisor.

First note that any $C\in \C L \cup \C B\cup \{Q\}$, with $C^2=-1$ is effective.  So, write
$C = f_1F_1 + \cdots + f_sF_s$ with the $F_i$ effective and prime and the $f_i$ nonnegative integers. Now
$$
-1=C^2= f_1C\cdot F_1 + \cdots + f_sC\cdot F_s
$$
and so $C \cdot F_i < 0$ for some $i$ with $f_i>0$.  But
$$
C\cdot F_i = f_1F_1\cdot F_i + \cdots + f_iF_i^2 + \cdots + f_sF_s\cdot F_i
$$
and $F_j\cdot F_i \geq 0$ for all $j \neq i$.  Thus $F_i^2 < 0$ as well.

Now suppose that $C\cdot D \geq 0$ for all $D \in \nneg(X)$.  If $F_i^2 \leq -2$ then $F_i \in \nneg(X)$  and so $C\cdot F_i \geq 0$ which is a contradiction. So, we must have $F_i^2 = -1$.  Hence, by $(b)$,  $F_i\in  \C L \cup \C B\cup \{Q\}$. 
But we can write down every element of $\C L \cup \C B\cup \{Q\}$
with self-intersection $-1$ and explicitly check that no two of them meet negatively.
Thus $C\cdot F_i<0$ implies $C=mF_i$ for some positive multiple $m$ of $F_i$, but by explicit check,
if $mF_i$ is in $\C L \cup \C B\cup \{Q\}$, then $m=1$, and hence $C=F_i$ is a prime divisor.

As for part $(e)$, we already observed that $A_r$ meets every element in $  \C L \cup \C B\cup \{Q\}$ positively, and from $(b)$ we have $\Neg(X) \subset  \C L \cup \C B\cup \{Q\}$.  Since $\Neg(X)$ generates EFF$(X)$, that is enough to finish the proof.
\end{proof}

We now show how to use Proposition \ref{ptsonconic} to tell if $F$ is effective or not, and to calculate $h^0(X,F)$ when $F$ is effective --- {\it assuming that we know the finite set $\nneg(X)$ and hence $\Neg(X)$}.
 
Let $H = F,\ N = 0$.  If $H\cdot C < 0$ for some $C \in \Neg(X)$, replace $H$ by $H-C$
(note that this reduces $H\cdot A_r$ since $C\cdot A_r > 0$) and replace $N$ by $N = N+C$.

Eventually either $H\cdot A_r < 0$ (hence $H$ and thus $F$ is not effective) or
$H\cdot C \geq 0$ for all $C \in \Neg(X)$, hence $H$ is nef  and effective by Proposition \ref{ptsonconic} $(a)$,
  and we have a Zariski decomposition $F = H + N$ with $h^0(X,F) = h^0(X,H)$ and so
$h^0(X,F) = (H^2-H.K_X)/2 + 1$ by Theorem \ref{BettisEtc}.

\medskip  What becomes clear from this procedure is that if we start with two sets of $r$ points on a conic and the isomorphism $\phi$ of groups given in Definition \ref{negcurtype} takes $\nneg(X)$ bijectively to $\nneg(X^\prime)$ then the calculations of $h^0(X,F)$ and $h^0(X^\prime , \phi(F))$ are formally the same and give $h^0(X,F)= h^0(X^\prime , F)$. 

From the definition we see that the classes in $\nneg(X)$ are simply the classes of self-intersection $< -1$ in $\Neg(X)$, so $\Neg(X)$ determines $\nneg(X)$.  Part $(c)$ of Proposition \ref{ptsonconic} tells us that $\nneg(X)$ determines $\Neg(X)$.  Thus, it is clear that if we have two configurations of $r$ points on a conic, then $\phi$ takes $\Neg(X)$ bijectively onto $\Neg(X^\prime)$ if and only if it takes $\nneg(X)$ bijectively onto $\nneg(X^\prime)$.

It follows that two configurations of $r$ points on a conic have the same configuration type if and only if they have the same negative curve type.  That type is also completely determined by the subset of  $ \C L \cup \C B\cup \{Q\}$ that turns out to be $\nneg(X)$.

Thus, to enumerate all configuration types of $r$ points on a conic, it is enough to enumerate all the subsets of $ \C L \cup \C B\cup \{Q\}$ which can be $\nneg(X)$.

\medskip This classification turns out to be rather simple in the case of $r$ points on a conic.  It is this fact that distinguishes this case from that of $6, 7$ or $8$ points (which we will treat in the following sections).

We now work out the classification, up to labeling of the points. We distinguish between two possibilities for $r$ points on a conic: either the $r$ points lie on an irreducible conic, or they lie on a reducible conic and no irreducible conic passes through the points.

If the $r$ points are on an irreducible conic, then since no line meets an irreducible conic in more than 2 points, no elements of $\C L$ are ever in $\nneg(X)$.  Elements of $\C B$ are never in $\nneg(X)$, so we need only decide if $Q$ is in $\nneg(X)$. Thus there are two cases:
\begin{itemize}
\item[{\bf I.}]\ $\nneg(X)$ empty;  i.e. $Q\not\in\nneg(X)$.  In this case we must have $r \leq 5$.
\item[{\bf II.}]\ $\nneg(X) = \{ Q \}$; in this case $r > 5$.
\end{itemize}

Now assume that the $r$ points lie only on a reducible conic.  
In this case, $Q$ is never a prime divisor and so we must decide only which of the elements of $\C L$ are in $\nneg(X)$.
Note also that in this case we have $r \ge 3$ (and if $r=3$, the 3 points must be collinear).

So, suppose that the $r$ points consist of $a$ on one line and $b$ on another line (where we may assume that $0\leq a \leq b$).  
We have to keep in mind the possibility that one of the points could be on both lines, i.e. it could be the point of intersection of the two lines.  Thus $a$ and $b$ satisfy the formula $a + b = r + \epsilon$, where $\epsilon = 0$ if the point of intersection of the two lines is not among our $r$ points and 
$\epsilon = 1$ if it is. We thus obtain two additional cases (depending on whether or not $a\ge3$), 
with two subcases (depending on whether $\epsilon$ is 0 or 1) when $a\ge3$:
\begin{itemize}
\item[{\bf III.}] $a \geq 3$ with $0\le\epsilon\le1$: in this case $\nneg(X)$ contains exactly two classes. Up to relabeling, these
classes are $L  - E_{1} - \cdots - E_{a}$ and  $L- E_{a+1} - \cdots - E_{r}$ if $\epsilon=0$, and
$L  - E_{1} - \cdots - E_{a}$ and  $L- E_{1} - E_{a+1}-\cdots - E_{r}$ if $\epsilon=1$.
\item[{\bf IV.}] $a < 3,\ b \geq 3$ with $\epsilon=0$: 
in this case, up to relabeling, $\nneg(X)$ contains exactly one class, 
namely $L - E_{1} - \cdots - E_{b}$. (Note that 
$a=2$ with $\epsilon=1$ gives the same configuration type
as does $a=1$ with $\epsilon=0$, and that $a=1$ with $\epsilon=1$ gives the same configuration type
as does $a=0$ with $\epsilon=0$. Thus when $a<3$ we may assume $\epsilon=0$.)
\end{itemize}

\begin{rem}\label{possibilities} \rm
It is worth observing that here there were very few possibilities for $\nneg(X)$. The fact is that by Bezout considerations, certain prime divisors cannot `{\it coexist}' in $\nneg(X)$.  For example, since a line cannot intersect an irreducible conic in more than two points we couldn't have $Q$ and $L-E_{i_1} - E_{i_2} - E_{i_3}$ both in $\nneg(X)$ at the same time.  The fact that it is easy to decide what can and cannot coexist is what makes the case of points on a conic relatively simple.  In the next section we will see that the analysis of which divisors can coexist in $\nneg(X)$ is much more subtle and even depends on the characteristic of our algebraically closed field.
\end{rem}

There are thus only finitely many
configuration types for $r$ points on a conic. Given any such points $p_i$,
we also see by Theorem \ref{BettisEtc} and Proposition \ref{ptsonconic} that the coefficients
$m_i$ and the configuration type of the points
completely determine $h_Z$ and the graded Betti numbers
of $I(Z)$ for $Z=m_1p_1+\cdots+m_rp_r$.

\begin{example}\label{ptsonconicex} \rm
We now work out a specific example, showing how Bezout considerations
completely determine the Hilbert  function and graded Betti numbers
in the case of any $Z$ with support in a conic. (By \cite{GuH}, the same considerations
apply for any $Z$ with support at any $r\le 6$ points. If $7\le r\le8$,
Bezout considerations still determine the Hilbert function, but the Betti
numbers are currently not always known.)

By Bezout considerations we mean this: To determine Hilbert functions and
graded Betti numbers it is enough to determine fixed
components, and to determine fixed components, it is enough to 
apply Bezout's theorem (which says that two curves whose intersection
is more than the product of their degrees must have a common component).

Consider distinct points $p_1,\ldots,p_5$ on a line $L_1$,
and distinct points $p_6,\ldots,p_9$ on a different line $L_2$ such that none of the points
is the point $L_1\cap L_2$ which we take to be $p_{10}$.
Let $Z=4p_1+2(p_2+\cdots+p_5)+3p_6+3p_7+2p_8+2p_9+3p_{10}$.
Then $Z$ has degree 46 and $h_Z(t)$ for $0\le t\le 15$ is
1, 3, 6, 10, 15, 21, 28, 31, 35, 38, 40, 42, 44, 45, 46, 46.
The complete set of graded Betti numbers for $F_0$ is
$7^5$, $9^3$, $10^1$, $13^1$, $15^1$, and for $F_1$ they are
$8^5$, $9^1$, $10^1$, $11^1$, $14^1$, $16^1$, where,
for example, $7^5$ for $F_0$ signifies that $t_7=5$, and
$8^5$ for $F_1$ signifies that $s_8=5$.

We first demonstrate how to compute $h_Z(t)$. 
First, $\binom{t+2}{2} - h_Z(t)=h^0(X, F(Z, t))$, where $X$ is obtained by blowing up the points
$p_i$, and $F(Z,t)=tL-4E_1-2(E_2+\cdots+E_5)-3E_6-3E_7-2E_8-2E_9-3E_{10}$. 
Clearly the proper transform $L_1'=L-E_1-\cdots-E_5-E_{10}$ 
of $L_1$ is in $\Neg(X)$, as are $L_2'=L-E_6-\cdots-E_9-E_{10}$,
$L_{ij}'=L-E_i-E_j$ for  $1\le i\le 5$ and $6\le j\le 9$, and $E_i$ for $1\le i\le 10$. 
Since the proper transform of no other line has negative self-intersection, we
see by Proposition \ref{ptsonconic} that $\Neg(X)$ consists exactly of $L_1'$, $L_2'$, the $L_{ij}'$
and the $E_i$.

Consider $h(6)$. Then $F(X,6)\cdot L_1'=-9$; i.e., any curve of degree
6 containing $Z$ has an excess intersection with $L_1$ of 9:
we expect a curve of degree 6 to meet a line only 6 times,
but containing $Z$ forces an intersection of 4 at $p_1$, 2 at $p_2$, etc,
for a total of 15, giving an excess of 9. Thus any such curve must contain
$L_1$ as a component. Said differently, if $F(Z,6)$ is effective, then so
is $F(Z,6)-L_1'=5L-3E_1-1(E_2+\cdots+E_5)-3E_6-3E_7-2E_8-2E_9-2E_{10}$.
But $(F(Z,6)-L_1')\cdot L_2'<0$, hence if $F(Z,6)-L_1'$ is effective , so is
$F(Z,6)-L_1'-L_2'$. Continuing in this way, we eventually conclude
that if $F(Z,6)$ is effective, so is $F(Z,6)-3L_1'-2L_2'-(E_2+\cdots+E_5)-2E_{10}-L_{1,6}=
-E_7$. But $A_{10}\cdot (-E_7)<0$, 
so we know $F(Z,6)-(3L_1'+2L_2'+(E_2+\cdots+E_5)+2E_{10}-L_{1,6})$ is not effective,
so $h_Z(6)=28$. 

As another example, we compute $h_Z(12)$. 
As before, if $h^0(X, F(Z,12))=h^0(X, H)$ where
$H=F(Z,12)-L_1'-L_2'=
10L-3E_1-1(E_2+\cdots+E_5)-2E_6-2E_7-E_8-E_9-E_{10}$. 
But $H$ is nef, since $H$ meets every element of
$\Neg(X)$ nonnegatively, so
$h^0(X, H)=\binom{10+2}{2}-19=47$ by Theorem \ref{BettisEtc}(a).
Thus $h_Z(12)=91-47=44$. Also, the linear system
$F(Z,12)$ decomposes as $H+(L_1'+L_2')$, 
where $H$ is nef and $L_1'+L_2'$ is fixed.
A similar calculation shows $F(Z,13)$ has fixed part
$L_1'$. Since the complete linear system of curves of degree 12 through $Z$
has a degree 2 base locus (corresponding to $L_1'+L_2'$), 
while the complete linear system of curves of degree 13 through $Z$
has only a degree 1 base locus (i.e., $L_1'$), the map
$\mu_{12}:I(Z)_{12}\otimes R_1\to I(Z)_{13}$
cannot be surjective, so we know at least 1 homogeneous generator
of $I(Z)$ is required in degree 13. But $\mu_{H}$ is surjective,
by Theorem \ref{BettisEtc}. Thus $\dim\cok(\mu_{12})
= h^0(X, F(Z,13))-h^0(X, H+L)=60-59=1$. Thus the graded Betti number for $F_0$ in
degree 13 is exactly 1. We compute the rest of the graded Betti numbers for
$F_0$ the same way. The graded Betti numbers for $F_1$
can be found using the well known formula
$s_i=t_i+\Delta^3h_Z(i)$ for $i>0$ and the fact that $s_0=0$.
\end{example}

\begin{rem}\label{ptsonconicrem} \rm 
One problem raised in \cite{GMS} is the existence and determination of maximal and minimal
Hilbert functions. For example, \cite{GMS} shows that there must be some 
$Z'$ such that $h_{2Z'}$ is at least as big in every degree as $h_{2Z}$ for
every $Z$ with $h_Z=h_{Z'}$; this $h_{2Z'}$ is referred to as ${\underline h}^{max}$.
The proof in \cite{GMS} is nonconstructive, and \cite{GMS}
determines ${\underline h}^{max}$ in only a few special cases.
The paper \cite{GMS} also raises the question of whether 
${\underline h}^{min}$ always exists; i.e., whether there exists a $Z'$
such that $h_{2Z}$ is at least as big in every degree as $h_{2Z'}$ for
every $Z$ with $h_Z=h_{Z'}$. This question remains open.

Below, we will consider all possible Hilbert functions $h_Z$, where $Z=p_1+\cdots+p_r$ for
points $p_i$ on a conic, and for each possible Hilbert function $h$, we will determine which 
configuration types $Z$ have $h_Z=h$. We will also determine $h_{2Z}$ for each
configuration type $Z$. One could, as shown in Example \ref{ptsonconicex},
also determine the graded Betti numbers. Thus, for points on a conic, this completes the program
begun in \cite{GMS} of determining the Hilbert functions (and graded Betti numbers) of double point schemes
whose support has given Hilbert function. Moreover,
by an inspection of the results below, one sees that among the configuration types 
of reduced schemes $Z=p_1+\cdots+p_r$ for points on a conic with a given Hilbert function
there is always one type which specializes to each of the others, and there is always one which is a specialization of each of the others. It follows by semincontinuity that $\underline{h}^{max}$ and $\underline{h}^{min}$ exist,
as do the analogous functions {\it for any m}; i.e., not only for the schemes $2Z$ but for the schemes $mZ$ for every $m\ge1$. Not only do these minimal and maximal Hilbert functions exist,
but because there are only finitely many configuration types for points on a conic
they can be found explicitly for any specific $m$: in particular, pick any
subscheme $Z=p_1+\cdots+p_r$ of distinct points $p_i$ such that $h_Z(2)<6$.
Since $h_Z(2)<6$, for any $Z'=p_1'+\cdots+p_r'$ such that $h_{Z'}=h_Z$, the points
$p_i'$ must lie on a conic. Now, for any given $m$, compute $h_{mZ'}$ (as demonstrated in
Example \ref{ptsonconicex}) for each configuration type and choose
those $Z'$ for which $h_{2Z'}$ is maximal, or minimal, as desired, from among
all $Z'$ for which $h_{Z'}=h_Z$. In addition to finding the maximal and minimal Hilbert functions,
this method also determines which configuration types give rise to them.

Now we determine the configuration types corresponding to each
Hilbert function $h$ for reduced 0-dimensional subschemes
with $h(2)<6$ (i.e., $h=h_Z$ for subschemes 
$Z=p_1+\cdots+p_r$ for distinct points $p_i$ on a conic),
following which we give the Hilbert function $h_{2Z}$ for each type $Z$.
We do not give detailed proofs; instead we note that
in case the points $p_i$ are smooth points on a reducible conic,
the Hilbert functions of $Z$ and of $2Z$ can be easily written down using 
the results of \cite{GMR} for $Z$ or using \cite{GMS} (for $\hbox{char}(k)=0$)
or \cite{CHT} (for $\hbox{char}(k)$ arbitrary) for either $Z$ or $2Z$. 
With somewhat more effort, the case of points on an irreducible conic
and the case of points on a reducible conic when one of the points $p_i$
is the singular point of the conic can be analyzed using the method
of Example \ref{ptsonconicex}.

It is convenient to specify a Hilbert function
$h$ using its first difference, which we will denote $\Delta h$.
Using this notation, the possible Hilbert functions for $Z=p_1+\cdots+p_r$
for distinct points $p_i$ on a conic are precisely those of the form
$$\Delta h=1\ \underbrace{2\ \cdots \ 2}_{i \hbox{\ times}}\ 
\underbrace{1\ \cdots \ 1}_{j \hbox{\ times}}\ 0\ \cdots,$$
where $i\ge0$ and $j\ge0$. It turns out that every Hilbert 
function $h$ of $r$ points on a conic of $\pr2$ corresponds to 
either one, two, three or four configuration types of $r$ points.  
We now consider each of the possibilities.

We begin with the case that $i=0$ and $j\ge0$:
$$\Delta h=1\ \underbrace{1\ \cdots \ 1}_{j}\ 0\ \cdots\ .$$
Here $h(1)<3$, so in this case the $r$ points are collinear (with $r=j+1$), hence there is
only one configuration type, it has $\nneg(X)=\emptyset$ for $r<3$
and $\nneg(X)=\{L - E_{1} - \cdots - E_{r}\}$ for $r\ge3$,
and $\underline{h}^{min} = \underline{h}^{max}$.

For the remaining cases we have $i>0$, so $h(1)=3$ but $h(2)<6$,
hence the $r$ points $p_i$ lie on a conic, but not on a line.
We begin with $i=1$ and $j>0$, so $r=j+3$:
$$
\Delta h = 1 \ 2\ \underbrace{1\ \cdots\ 1}_{j}\  0\ \cdots\ .
$$
We obtain this Hilbert function either for $r=4$ general points (in which case $\nneg(X)=\emptyset$),
or in case we have $r-1$ points on one line and one point off that line (in which case, up to equivalence, 
$\nneg(X)=\{L - E_{1} - \cdots - E_{r-1}\}$). 
Thus, when $r=4$ this Hilbert function arises from two configuration types, and for $r>4$
it arises from exactly one configuration type. Thus except when $r=4$, we have $\underline{h}^{min} = \underline{h}^{max}$, and, when $r=4$, $\underline{h}^{max}$ occurs when the four points are general while $\underline{h}^{min}$ occurs
when exactly three of the four are collinear.

Next consider $r = 2i+j+1$, $i\geq 2$, $j \geq 2$.
$$
\Delta h =  1\ \underbrace{2\ \cdots\ 2}_{i}\ \underbrace{1\ \cdots\ 1}_{j}\ 0\ \cdots\ .
$$
This can occur for points only on a reducible conic; in this case, recalling the formula $a+b=r+\epsilon$, 
we have $b=i+j+1$ and $a=i+\epsilon$. 
If $i=2$ and $\epsilon=0$, then 
$\nneg(X)=\{L - E_{i_1} - \cdots - E_{i_b}\}$. Otherwise, we have
$\nneg(X)=\{L - E_{i_1} - \cdots - E_{i_b}, L - E_{i_{b+1}} - \cdots - E_{i_{b+a}}\}$, with
the indices being distinct if $\epsilon=0$ but with $i_1=i_{b+1}$ if $\epsilon=1$.
Thus this Hilbert function arises from two configuration types for each given $i$ and $j$,
and $\underline{h}^{min}$ occurs when $\epsilon=1$ and $\underline{h}^{max}$ when $\epsilon=0$.

Now consider $i>0$ and $j=0$, so $r = 2i+1$:
$$
\Delta h = 1\ \underbrace{2\ \cdots \ 2}_{i}\ 0\ \cdots\ .
$$
This case comes either from the configuration type with 
$r\ge3$ points (with $r$ odd) on an irreducible conic (in which case $\nneg(X)=\emptyset$ if $i=1$ or 2, and
$\nneg(X)=\{2L - E_{1} - \cdots - E_{r}\}$ if $i>2$), or
from the configuration types with $r=a+b-\epsilon$ points on a pair of lines,
for which $b=i+1$, $a=i>1$ and $\epsilon=0$
(in which case, up to equivalence, $\nneg(X)=\{L - E_{1} - \cdots - E_{3}\}$
if $i=2$ and $\nneg(X)=\{L - E_{1} - \cdots - E_{b}, L - E_{b+1} - \cdots - E_{r}\}$
if $i>2$) or for which
$b=a=i+1>2$ and $\epsilon=1$
(in which case $\nneg(X)=\{L - E_{1} - \cdots - E_{b}, 
L - E_1 - E_{b+1} - \cdots - E_{r}\}$). 
Thus this Hilbert function arises from a single configuration type if $i=1$
(in which case $\underline{h}^{min} = \underline{h}^{max}$),
and from three configuration types for each $i>1$
(in which case $\underline{h}^{min}$ occurs when the conic is reducible with
$b=a$, and $\underline{h}^{max}$ occurs when the conic is irreducible).

Finally, consider $i>1$ and $j=1$ (note that $i=j=1$ is subsumed by the previous case), so $r=2i+2$:
$$
\Delta h =  1\  \underbrace{2\ \cdots\ 2}_{i}\ 1\ 0\ \cdots\ .
$$
This case comes from either the configuration type with 
$r\ge6$ points (with $r$ even) on an irreducible conic 
(in which case $\nneg(X)=\{2L - E_{1} - \cdots - E_{r}\}$), or 
from the configuration types
with $r=a+b-\epsilon$ points on a pair of lines
for which either: $b=i+2$, $a=i$ and $\epsilon=0$
(in which case $\nneg(X)=\{L - E_{i_1} - \cdots - E_{i_b}\}$
if $i=2$ and $\nneg(X)=\{L - E_{i_1} - \cdots - E_{i_b}, L - E_{i_{b+1}} - \cdots - E_{i_{b+a}}\}$
if $i>2$); 
$b=a=i+1$ and $\epsilon=0$
(in which case $\nneg(X)=\{L - E_{i_1} - \cdots - E_{i_b}, L - E_{i_{b+1}} - \cdots - E_{i_{b+a}}\}$); or
$b=i+2$, $a=i+1$ and $\epsilon=1$
(in which case $\nneg(X)=\{L - E_{1}- E_{i_1}- \cdots - E_{i_b}, L - E_{1}- E_{i_{b+1}} - \cdots - E_{i_{b+a}}\}$).
Thus this Hilbert function arises from four configuration types for each $i$, and
$\underline{h}^{min}$ occurs when the conic is irreducible
while $\underline{h}^{max}$ occurs when the conic is reducible with $\epsilon=1$.

We now write down the Hilbert functions of double point schemes
for each configuration type. In the case of $r$ points on a line,
the corresponding double point scheme $2Z$
has Hilbert function
$$\Delta h_{2Z}= 1\ \underbrace{2\ \cdots \ 2}_{r}\ \underbrace{1 \ 1\ \cdots \ 1}_{r-1}\ 0\ \cdots.$$

The possible Hilbert functions for $2Z$ for $Z=p_1+\cdots+p_r$ contained in a conic
but not on a line are all of the form
$$\Delta h_{2Z}=1\ 2\ 3\ \underbrace{4\ \cdots\ 4}_i\ \underbrace{3\ \cdots\ 3}_j\ 
\underbrace{2\ \cdots\ 2}_k\ \underbrace{1\ \cdots\ 1}_l\ 0\ \cdots,$$
for various nonnegative integers $i,j,k$ and $l$.

First suppose the conic is irreducible. Since $r=2$ is the case of points on a line, we may assume
that $r\ge 3$. If $r=3$, then
$$\Delta h_{2Z}=1\ 2\ 3\ 3\ 0\ \cdots.$$
If $r=4$, then
$$\Delta h_{2Z}=1\ 2\ 3\ 4\ 2\ 0\ \cdots.$$
If $r\ge 5$ is odd, then 
$$\Delta h_{2Z}=1\ 2\ 3\ \underbrace{4\ \cdots\ 4}_{i=\frac{r-1}{2}}\ 
\underbrace{2\ \cdots\ 2}_{k=\frac{r-5}{2}}\ 1\ 0\ \cdots.$$
If $r\ge 6$ is even, then 
$$\Delta h_{2Z}=1\ 2\ 3\ \underbrace{4\ \cdots\ 4}_{i=\frac{r}{2}-1}\ 
3\ \underbrace{2\ \cdots\ 2}_{k=\frac{r-6}{2}}\ 1\ 0\ \cdots.$$

Now assume the points lie on a reducible conic, consisting of two distinct lines, and that 
no line nor any irreducible conic contains the points.
Thus $4\le r=a+b-\epsilon$, and we may assume that $1\le a\le b$.

If $2a<b$ and $\epsilon=0$, then
$$\Delta h_{2Z}=1\ 2\ 3\ 
\underbrace{4\ \cdots\ 4}_{i=a}\ 
\underbrace{3\ \cdots\ 3}_{j=a-1}\ 
\underbrace{2\ \cdots\ 2}_{k=b-2a-1}\ 
\underbrace{1\ \cdots\ 1}_{l=b-1}\ 0\ \cdots.$$
If $a<b<2a$ and $\epsilon=0$, then
$$\Delta h_{2Z}=1\ 2\ 3\ 
\underbrace{4\ \cdots\ 4}_{i=a}\ 
\underbrace{3\ \cdots\ 3}_{j=b-a-1}\ 
\underbrace{2\ \cdots\ 2}_{k=2a-b-1}\ 
\underbrace{1\ \cdots\ 1}_{l=2b-2a-1}\ 0\ \cdots.$$
If $b=2a$ and $\epsilon=0$, then
$$\Delta h_{2Z}=1\ 2\ 3\ 
\underbrace{4\ \cdots\ 4}_{i=a}\ 
\underbrace{3\ \cdots\ 3}_{j=a-2}\ 2\ 
\underbrace{1\ \cdots\ 1}_{l=2(a-1)}\ 0\ \cdots.$$
The only remaining case with $\epsilon=0$ is $b=a\ge3$, in which case
$$\Delta h_{2Z}=1\ 2\ 3\ \underbrace{4\ \cdots\ 4}_{i=a-1}\ 
3\ \underbrace{2\ \cdots\ 2}_{k=a-3}\ 1\ 0\ \cdots.$$

If $\epsilon=1$, we may as well assume that $3\le a\le b$,
since $a=2$ and $\epsilon=1$ is the same as $a=1$
and $\epsilon=0$ (for appropriate $b$'s).
If $a=b$ and $\epsilon=1$, then
$$\Delta h_{2Z}=1\ 2\ 3\ 
\underbrace{4\ \cdots\ 4}_{i=a-2}\ 3\ 
\underbrace{2\ \cdots\ 2}_{k=a-2}\ 0\ \cdots.$$
If $a<b\le2a-1$ and $\epsilon=1$, then
$$\Delta h_{2Z}=1\ 2\ 3\ 
\underbrace{4\ \cdots\ 4}_{i=a-1}\ 
\underbrace{3\ \cdots\ 3}_{j=b-a-1}\ 
\underbrace{2\ \cdots\ 2}_{k=2a-b-1}\ 
\underbrace{1\ \cdots\ 1}_{l=2(b-a)}\ 0\ \cdots.$$
If $2a-1\le b$ and $\epsilon=1$, then
$$\Delta h_{2Z}=1\ 2\ 3\ 
\underbrace{4\ \cdots\ 4}_{i=a-1}\ 
\underbrace{3\ \cdots\ 3}_{j=a-2}\ 
\underbrace{2\ \cdots\ 2}_{k=b-2a+1}\ 
\underbrace{1\ \cdots\ 1}_{l=b-1}\ 0\ \cdots.$$
\end{rem}

\section{The Case of $r \le 8$ distinct points of $\pr2$}\label{7and8pts}

We now consider the case of $r\le8$ distinct points of $\pr2$.  As in the case of points on a conic, we define some finite families of special divisors on the surface $X$ obtained by blowing up the points $p_1, \ldots , p_r$, where $r \leq 8$.  As in earlier sections we will let $E_i \in \hbox{Cl}(X)$ be the divisors corresponding to the blown up points and $L$ the divisor class which is the preimage of a general line in $\pr2$.  The families will be denoted: 
\begin{itemize}
  \item ${\mathcal B}_r = \{ E_1, \ldots , E_r \}$;
  \item ${\mathcal L}_r = \{ L - E_{i_1}-\cdots -E_{i_j} \ \mid \ 2 \leq j \leq r \}$;
  \item ${\mathcal Q}_r = \{ 2L - E_{i_1} - \cdots - E_{i_j} \ \mid \ 5 \leq j \leq r \}$;
  \item ${\mathcal C}_r = \{3L - 2E_{i_1} - E_{i_2} - \cdots - E_{i_j} \mid \ 7 \leq j \leq 8, j \leq r \}$; and
  \item ${\mathcal M}_8 =  \{ 4L - 2E_{i_1}-2E_{i_2}-2E_{i_3}-E_{i_4} - \cdots - E_{i_8}, 5L - 2E_{i_1} - \cdots - 2E_{i_6}-E_{i_7}-E_{i_8}, \newline6L-3E_{i_1}-2E_{i_2} - \cdots - 2E_{i_8} \}$.
\end{itemize}

Let $\C N_r = \C B_r \cup \C L_r\cup \C Q_r \cup \C C_r \cup \C M_8$.
Notice that the elements of ${\mathcal B}_r$ and ${\mathcal M}_8$
all have self-intersection $-1$, hence are never elements of
$\nneg(X)$. They will be elements of $\Neg(X)$ when they are classes of
irreducible curves and as such are involved in determining Hilbert functions
via the role they play in determining the fixed and free parts
of linear systems $|F(Z,t)|$. (See Example \ref{ptsonconicex}
for an explicit demonstration of how $\Neg(X)$ is involved 
in computing Hilbert functions.)

\begin{prop}\label{8pts} 
Let $X$ be obtained by blowing up $2\le r\le8$ distinct points of $\pr2$.
\begin{itemize}
\item[(a)] Then $\NEF(X)\subset \EFF(X)$;
\item[(b)] If $r=8$, then every element of $\EFF(X)$ is a nonnegative rational
linear combination of elements of $\Neg(X)$;
\item[(c)] If $r<8$, then every element of $\EFF(X)$ is a nonnegative integer
linear combination of elements of $\Neg(X)$;
\item[(d)] $\Neg(X)\subset \C N_r$;
\item[(e)] $\Neg(X)=\nneg(X) \cup \{C\in \C N_r\,\,\mid \,\,C^2 = -1, C\cdot D\ge 0\hbox{ for all }D\in \nneg(X)\}$;
\item[(f)] for every nef divisor $F$, $|F|$ is fixed component free; and
\item[(g)] the divisor $A_r$ is ample.
\end{itemize}
\end{prop}

\begin{proof}
Part (a) is Theorem 8 of \cite{BHProc}. Part (f) follows from \cite{ARS}.
For (d), apply adjunction: if $C$ is the class of a prime divisor with $C^2<0$,
then $C^2=-C\cdot K_X+2g-2$, where $g$ is the arithmetic genus of $C$.
Since $g\ge 0$, we see that either $-C\cdot K_X<0$ (and hence $C$ is a 
fixed component of $|-K_X|$), or $-C\cdot K_X\ge 0$ and hence 
$g=0$ (so $C$ is smooth and $g$ is actually the genus of $C$)
and either $-C\cdot K_X=-1=C^2$ or $-C\cdot K_X=0$ and $C^2=-2$.
But it is easy to check that $K_X^\perp$ is negative definite so it is not hard to
list every class satisfying $-C\cdot K_X=0$, $C^2=-2$ and $L\cdot C\ge0$.
If one does so, one gets the subset of $\C N_r$ of elements of self-intersection $-2$; 
see Chapter IV of \cite{Manin} for details.
Similarly, the set of all solutions to $C^2=-1=-C\cdot K_X$ consists of the elements in $\C N_r$ of self intersection $-1$;
again, we refer the reader to Chapter IV of \cite{Manin}.
If $-C\cdot K_X<0$, then $-K_X-C$ is effective so $C\cdot L\le -K_X\cdot L=3$.
It is easy to check that the only such curves $C=dL-m_1E_1-\cdots-m_rE_r$
with $d\le 3$ which are effective, reduced and irreducible,
come from lines through 4 or more points (i.e., up to indexation,
$C=L-E_1-\cdots-E_j$, $j\ge4$), or conics through 7 or 8 points
(i.e., $C=2L-E_1-\cdots-E_j$, $7\le j\le8$). There can be no cubic,
since a reduced irreducible cubic has at most one singular point, so 
would be either $-K_X$ itself or $3L-2E_1-E_2-\cdots-E_8$, neither of which
meets $-K_X$ negatively. Thus if $-C\cdot K_X<0$, then $C$
must be among the elements of $\C N_r$ with self-intersection less than $-2$.
This proves (d). The proof of (e) is essentially the same as that 
of Proposition \ref{ptsonconic}(c). 
Parts (b) and (c) are well known; see \cite{GuH} for a detailed
proof in case $r=6$. The proofs for other values of $r$ are similar. 
Finally, the proof of (g) is essentially the 
same as the proof of Proposition \ref{ptsonconic}(e). 
\end{proof}

With Proposition \ref{8pts}, the procedure for computing
$h^0(X, F)$ for any class $F$ on a blow-up $X$ of
$\pr2$ at $r\le 8$ points is identical to the procedure if the points lie on a conic.
Given $\nneg(X)$, find $\Neg(X)$. Let $H=F$. Loop through the elements $C$ of 
$\Neg(X)$. Whenever $F\cdot C<0$, replace $H$ by $H-C$.
Eventually either $H\cdot A_r<0$, in which case $h^0(X, F)=h^0(X, H)=0$,
or $H$ is nef and
$F=H+N$, where $N$ is the sum of the classes subtracted off,
and we have $h^0(X, F)=h^0(X, H)=(H^2-K_X\cdot H)/2+1$
by Theorem \ref{BettisEtc}.
In case $r\le 6$, we can also compute the graded Betti numbers,
as discussed above, since when $H$ is nef, $\mu_H$ has maximal
rank by Theorem \ref{BettisEtc}.

It follows, in the same way as for points on a conic, that
the negative curve type of sets of $r\le 8$ points is the same as
the configuration type.

\section{Configuration Types for $r=7$ or 8 Points of $\pr2$}\label{conftypes}

Given $r$ distinct points of $\pr2$, if either
$r\le 8$ or the points lie on a conic, then in previous sections we saw how to
compute the Hilbert function, and in some cases the graded Betti numbers,
of any fat point subscheme with support at the $r$ points, if we know
the configuration type of the points (or, equivalently, if we know
$\nneg(X)$, where $X$ is the surface obtained by blowing up the points).
The types were easy to enumerate
in the case of points on a conic. We now consider the case of $r\le 8$
points. (If $r>8$, then typically $\Neg(X)$ is infinite, and
$\nneg(X)$ can also be infinite. In fact, for each $r\ge9$, there are
infinitely many configuration types of $r$ points;
see Example \ref{infmanytypes}. Moreover, for $r>9$,
$\Neg(X)$ is known only in special cases, such as, for example, when 
$X$ is obtained by blowing up points on a conic.)
Since any $r\le 5$ points lie on a conic, and since the case $r=6$
is done in detail in \cite{GuH}, we will focus on the cases $7\le r\le 8$.

We begin by formalizing the notion of negative curve type: by Proposition \ref{8pts},
the elements of $\nneg(X)$ are all in $\C N_r$, and moreover if $C$ and $D$ are distinct
elements of $\nneg(X)$, then, being prime divisors, $C\cdot D\ge0$. 
Thus we define a subset $S$ of $\C N_r$ 
to be {\it pairwise nonnegative} if whenever $C\ne D$ for elements $C$ and $D$ of $S$,
then $C\cdot D\ge 0$. 

\begin{defn}\rm A {\it formal configuration type\/} for $r=7$ or 8 points
in $\pr2$ is a pairwise nonnegative subset $S$ of $\C N_r$.
We say two types $S_1$ and $S_2$ are equivalent if by 
permuting $E_1,\ldots,E_r$ we can transform $S_1$ to $S_2$.
If $S=\nneg(X)$ for some surface $X$, we say $S$ is 
{\it representable}; i.e., there is a bijective correspondence between
configuration types and representable formal configuration types.
\end{defn}

The formal configuration types for $r\le 5$ points are a special case of the types
for points on a conic, which are easy to enumerate and was done in Section \ref{PtsConic}.
In this section we will explicitly list the formal configuration types $S$ for sets of 6, 7 and 8 points,
determine which are representable, 
and for each such formal configuration type, represented say by
points $p_1,\ldots,p_r$, we will answer some of the questions raised in \cite{GMS}
by applying our procedure to compute Hilbert functions
(and, for $r\le 6$, the graded Betti numbers too),
for subschemes $Z=p_1+\cdots+p_r$ and $2Z$. Also, not only do we
find all the Hilbert functions that arise, we also determine exactly
what arrangements of the points are needed to give each Hilbert function.
We note that the same procedure 
can be used to compute the Hilbert function for $m_1p_1+\cdots+m_rp_r$
for any $m_i$, and what arrangements of the points give that Hilbert function,
as long as $r\le8$. 

\begin{rem}\label{combgeomrem} \rm Formal configuration 
types generalize the notion \cite{CR} of combinatorial geometries of rank up to 3.
A combinatorial geometry is a formal specification
of the linear dependencies on a set of ``points.''
Formally, a combinatorial geometry on $r$ ``points'' of rank at most 3
(hence formally the points are in the plane) can be defined
to be a matrix with $r$ columns. Each row contains only zeroes and ones,
the sum of the entries in each row is at least 3, and the dot product of two different
rows is never bigger than 1. The columns represent points, and each row represents
a line. A 1 in row $i$ and column $j$ means that line $i$ goes through 
point $j$; a 0 means that it doesn't. We allow a matrix with no rows, which just means
the given combinatorial geometry consists of $r$ points, no three of which are collinear.
However, writing down a matrix with the required properties doesn't
guarantee that any point set with the specified dependencies, and no other dependencies,
actually exists. Determining whether there exists an actual set of points exhibiting
a given combinatorial geometry is a separate question. When there
does exist one over a given field $k$, one says that the combinatorial geometry 
is {\it representable\/} over $k$. A list of all combinatorial geometries on 
$r\le 8$ points is given in \cite{BCH}, without dealing with the problem of
representability.

Note that if one takes a combinatorial geometry of rank at most 3 on $r\le8$ points, 
each row of its matrix gives an element of 
$\C N_r$ of self intersection $< -1$, as follows. Let $m_{ij}$ be the entry in column $j$ of
row $i$. Then $C_i=L-m_{i1}E_1-\cdots-m_{ir}E_r\in \C N_r$, and the set $S$
consisting of all of the $C_i$ is a formal configuration type. Conversely,
if $S$ is a formal configuration type such that every $C\in S$ has $C\cdot L=1$,
then we can reverse the process and obtain a 
combinatorial geometry from $S$. However, we are also interested
in cases where the points lie on conics and cubics, so $C\cdot L$ can, for
us, be 2 or 3. (We could also allow infinitely near points,
in which case we could have $C\cdot L=0$.)
\end{rem}

Since $S\subset \C N_r$ and $\C N_r$ is a finite set, it is clear
that with enough patience one can write down every
formal configuration type. The list given in \cite{BCH} gives us all
formal configuration types which do not involve curves of degree
bigger than 1. To these we add, in all possible ways, classes $C\in \C N_r$
with $C\cdot L=2$ or 3. Each time we find a class
that can be added to a previous formal configuration type
without violating pairwise nonnegativity, we obtain another type.
We must then check to see if it is equivalent to one which has already occurred.
Eventually we obtain a complete list. We now give these lists.

Table \ref{6pttableA} gives the list of the eleven formal configuration types 
(up to equivalence) for $r=6$ points.
It is not hard to see that each type is in fact representable.
In the table, the 6 points
$p_1,\ldots, p_6$ are denoted alphabetically by the lower case letters
``a'' through ``f''. Whenever a type $S$ contains a class
$L-E_{i_1}-\cdots-E_{i_j}$, the letters corresponding 
to the points $p_{i_1},\ldots,p_{i_j}$ are listed in the table,
preceded by a 1 to indicate that the coefficient of $L$ is 1.
Intuitively, configuration type 7, for example, is a specification that 
the points $p_1$, $p_2$, $p_3$ and $p_4$ are to be collinear, and so
are the points $p_1$, $p_5$ and $p_6$. For type 11, no three of the points
are to be collinear, but there is to be an irreducible conic containing the points
$p_1,\ldots,p_6$. Thus the type 11 lists ``2: abcdef'', and $S$ in this case consists
of $2L-E_1-\cdots-E_6$. Our notation mimics that of \cite{BCH}.
(Table \ref{6pttableA} also gives the types for $r<6$, using the following convention.
The types for $r=5$ correspond exactly to the types listed in 
Table \ref{6pttableA} which do not involve the letter ``f''; the types for $r=4$
are those which do not involve the letters ``e'' or ``f'', etc.)

Table \ref{6pttableB} gives both $h_Z$ and the graded Betti numbers for $Z=m(p_1+\cdots +p_r)$
for each configuration type with $r\le 6$ for $m=1$ and $m=2$. 
The configuration types are listed using the same numbers as in the previous table.
The table gives the Hilbert function $h_Z$ 
by listing $h_Z(t)$ for every $t$ from 0 up to the degree $t$ for which
$h_Z(t)=\hbox{deg}(Z)$.
The graded Betti numbers are specified using the same notation explained
in Example \ref{ptsonconicex}. 

Table \ref{7pttableA} lists the configuration types 
(up to equivalence) for $r=7$ points.
There are 29 types, and each is representable over some field, although as we shall see
in some cases representability depends on the characteristic of $k$.

Table \ref{7pttableB} gives the Hilbert functions for 
$Z=p_1+\cdots+p_7$ and $2Z$, for each type.
The configuration types are listed by
the same item number used in the table above.
Five different Hilbert functions occur for 7 distinct points of multiplicity 1.
For three of these five Hilbert functions $h_Z$, only one 
Hilbert function is possible for $h_{2Z}$.
For one of these five Hilbert functions, three different Hilbert 
functions occur for double points, and for the other,
two different Hilbert functions occur for double points. 
All together,
there are thus eight different Hilbert functions which occur
for 7 distinct points in the plane of multiplicity 2.
For each Hilbert function $h$ of 7 simple points, we see from the table
that there is both a maximum and a minimum Hilbert function among
all Hilbert functions of double points whose support has the given Hilbert function $h$.
The table groups together Hilbert functions of double point schemes whose
support schemes have the same Hilbert function.

Tables \ref{8pttableA1} and \ref{8pttableA2} list the 146 formal 
configuration types (up to equivalence) for $r=8$ points. 
Note that in every case of a configuration type involving a cubic,
the cubic has a double point which is always assumed to be at the last point;
i.e., the notation \hbox{3: abcdefgh} denotes a cubic through all 8 points,
always with a double point at h.

Table \ref{8pttableB} gives $h_Z$ and $h_{2Z}$ for each reduced scheme $Z$
corresponding to a formal configuration type of $r=8$ points. 
Six different Hilbert functions occur for 8 distinct points of multiplicity 1.
Again we see that there is both a maximum and a minimum Hilbert function
among the Hilbert functions of double point schemes 
whose support scheme's Hilbert function is one of these six.
For four of the six Hilbert functions $h_Z$, only one 
Hilbert function is possible for $h_{2Z}$.
For both of the other two of these six Hilbert functions, three different Hilbert 
functions occur for double points. 
All together,
there are thus ten different Hilbert functions which occur
for 8 distinct points in the plane of multiplicity 2.
(Note that types 30, 45 and 96 are not representable over any field $k$.
Our procedure can still be run, however, so
we show the result our procedure gives in these cases too.
This is what would happen, if these types were representable.)

\begin{example}\rm
We give an example to demonstrate how we generated the lists of configuration types. Consider the case of 8 points. From the list of the 69 simple eight point matroids of rank at most three given in \cite{BCH},
we can immediately write down all formal configuration types contained in $\C L_8$. This gives 69 configuration types.

For each of these 69 we then check to see which classes corresponding to either conics through 6 or more points or cubics through all 8 points and singular at one are compatible with the classes in the formal configuration type of the given matroid. 

For example, take configuration 10. Configuration 10 is contained in $\C L_8$ and thus comes directly from one of the 69 matroids. In terms of the basis $L,E_1,\cdots,E_8$, the classes 
in configuration type 10 are $L_i$ for $i=1,\ldots,4$
given as follows:
$$\begin{array}{rcrrrrrrrrrrrrrrrr}
 L_1: & = & 1 &-1& -1 &-1 & 0 & 0  &0  &0 & 0 \cr
L_2: & = & 1& -1 & 0 & 0& -1 &-1  &0 & 0 & 0 \cr
L_3: & = & 1 & 0 &-1 & 0 &-1 & 0 &-1 & 0 & 0 \cr
L_4: & = & 1 & 0 & 0 &-1 & 0 &-1&  0 &-1 & 0
\end{array}$$

After permuting the coefficients of the $E_i$,
there are 28 classes in the orbit of $2L-E_1-\cdots-E_6$,
eight of the form $2L-E_1-\cdots-E_7$, one of the form
$2L-E_1-\cdots-E_8$, and eight of the form
$3L-2E_1-E_2-\cdots-E_8$. Of these 45 classes, the only ones
which meet all four classes of configuration type 10 nonnegatively
are $Q_1$, $Q_2$ and $C$ given as follows:

$$\begin{array}{lcrrrrrrrrrrrrrrrrrr}
Q_1:& = & 2 &-1 & 0& -1 &-1 & 0& -1 &-1 &-1\cr
Q_2: &= & 2& -1& -1 & 0 & 0 &-1 &-1 &-1& -1\cr
C \ :& = & 3 &-1 &-1& -1 &-1& -1 &-1 &-1 &-2
\end{array}$$

Since $Q_1\cdot Q_2\ge0$ but $Q_i\cdot C<0$ for $i=1,2$,
the formal configuration types one gets from configuration 10
are: 
$\{L_1,L_2,L_3,L_4\}$ (i.e., type 10 itself),
$\{L_1,L_2,L_3,L_4, Q_1\}$ (this has type 77),
$\{L_1,L_2,L_3,L_4,Q_2\}$ (this also has type 77),
$\{L_1,L_2,L_3,L_4, Q_1, Q_2\}$ (this has type 89) and
$\{L_1,L_2,L_3,L_4, C\}$ (this has type 40).
Of course, none of these are isomorphic
to any obtained starting from a matroid
not isomorphic to the one giving configuration type 10,
but in fact $\{L_1,L_2,L_3,L_4, Q_1\}$ is isomorphic to
$\{L_1,L_2,L_3,L_4,Q_2\}$ (since permuting the basis
$L,E_1,\ldots,E_8$ by transposing indices 2 and 4 and
at the same time transposing indices 3 and 5 
permutes the $E_i$ in such a way
as to convert $\{L_1,L_2,L_3,L_4, Q_1\}$ into
$\{L_1,L_2,L_3,L_4,Q_2\}$).
\end{example}

We now consider the problem of representability. 
Note that every 7 point formal configuration type is also, clearly,
an 8 point formal configuration type. For example, type 24 for 7 points is
(after a permutation of the points; just swap d and g)
type 30 for 8 points. This type is also a combinatorial geometry:
it is the well known Fano plane; that is, the projective plane over 
the integers modulo 2. In particular, regarded as a 7 point formal configuration type,
it is representable in characteristic 2. Regarded as an 8 point type, 
however, it is never representable, since by the proof of Theorem \ref{8ptrepthm} no matter
where one picks the eighth point to be, it is either on one of the lines through
pairs of the first seven points or there is a singular cubic
through the first seven points with its singularity at the eighth
point. 

\begin{thm}\label{7ptrepthm} Let $k$ be an algebraically closed field.
Type  23 for $r=7$ points is representable
if and only if $\Char(k)\neq2$. 
Type  24 for $r=7$ points is representable
if and only if $\Char(k)=2$. 
The remaining types for $r=7$ are representable over
every algebraically closed field.
\end{thm}

Some lemmas will be helpful for the proof.
Let $S$ be a formal configuration type for $r$ points.
If $C=dL-m_1E_1-\cdots-m_rE_r\in S$ and $i\le r$,
let $C_i=dL-m_1E_1-\cdots-m_iE_i$ (i.e., $C_i$ is obtained from
$C$ by truncation). Let 
$S(i)$ be the set of all classes $C_i$ such that
$C\in S$ and $C_i^2\le-2$.

\begin{lem}\label{7ptreplem} Let $k$ be an algebraically closed field,
and let $S$ be a formal configuration type for $r=7$
points. If $S(i)$ is $k$-representable for some $i<r$,
and if the number $\#(S(i+1)-S(i))$ of elements 
in $S(i+1)$ but not in $S(i)$ is at most 1, 
then $S(i+1)$ is $k$-representable. In particular, 
if $\#(S(i+1)-S(i))\le 1$ for each $i<r$, then $S$ is representable
over $k$.
\end{lem}

\begin{proof}
Every element of $S(i)$ is the truncation of an element of $\C N_{i+1}$. 
Thus $S(1)$ is always empty, because there are no elements of $\C N_2$ of self intersection $< -1$.
Thus $S(1)$ is a $k$-representable 1 point formal configuration type.
Suppose $S(i)$ is a $k$-representable $i$ point formal configuration type.
Thus there are points $p_1,\ldots,p_i$ in the projective plane
$\pr2$ over $k$ such that for the surface $X_i$ obtained by 
blowing them up we have $S(i)=\nneg(X)$.

Now assume $S(i+1)=S(i)$.  Suppose for some point $p_{i+1}$, the surface $X_{i+1}$
obtained by blowing up $p_1,\ldots,p_{i+1}$ has a prime divisor
$C=dL-m_1E_1-\cdots-m_{i+1}E_{i+1}\in\nneg(X_{i+1})$ not in $S(i+1)$.
Then $C\in \C N_{i+1}$ with $C^2 < -1$.  But the coefficients $m_j$ are 
always 0 or 1 (since $i+1\le r=7$). Since $C\not\in S(i+1)=S(i)$, we see
$m_{i+1}=1$. Thus $C'=dL-m_1E_1-\cdots-m_{i}E_{i}$ 
is the class of a prime divisor on $X_i$ with $(C')^2\le-1$.
Thus $C'\in\Neg(X_i)$. Since $\Neg(X_i)$ is finite, 
the union of all such $C'$ is a proper closed subset of $X_i$.
If we pick $p_{i+1}$ to avoid this closed subset, then 
$\nneg(X_{i+1})=S(i+1)$, so $S(i+1)$ is representable over $k$.

Finally, assume $\#(S(i+1)-S(i))=1$. 
Let $D=dL-m_1E_1-\cdots-m_{i+1}E_{i+1}$ be the class
in $S(i+1)$ which is not in $S(i)$. Thus $m_{i+1}=1$
and $D'=dL-m_1E_1-\cdots-m_iE_i$ is in  $\C N_i$. If $(D')^2< -1$, then (keeping in mind
how $S(i)$ is constructed) $D'\in S(i)$, and hence $D'$ is the class of
a prime divisor on $X_i$. If $(D')^2 = -1$, then since
$D'\cdot C\ge0$ for all $C\in S(i)$, 
by Proposition \ref{8pts}(e) $D'$ again is the class of
a prime divisor on $X_i$. Thus we must choose $p_{i+1}$
to be a point on $D'$, which, as before, is not on any $C'$
in $\Neg(X_i)$ other than $D'$. These $C'$ comprise a proper
closed subset meeting $D'$ in a finite subset, so choosing
$p_{i+1}$ to be any point of $D'$ not in this finite subset
results in $\nneg(X_{i+1})=S(i+1)$, 
showing that $S(i+1)$ is $k$-representable.
\end{proof}

\begin{lem}\label{Torsionlem} Let $k$ be an algebraically closed field,
and let $S$ be a formal configuration type for $r\le8$
points. Assume that $S\subset K^\perp$, where $K=-3L+E_1+\cdots+E_r$
(i.e., $D\cdot K=0$ for every element $D\in S$).
Let $T$ be the torsion subgroup of $K^\perp/\langle S\rangle$, where
$\langle S\rangle$ is the subgroup of $K^\perp$ generated by $S$.
If for some reduced irreducible plane cubic $C\subset\pr2(k)$
either $K^\perp/\langle S\rangle=T$ and $T$ is isomorphic to a subgroup of 
$\Pic^0(C)$ or
there are infinitely many positive integers $l$ such that
$T\times {\bf Z}/l{\bf Z}$ is isomorphic to a subgroup of 
$\Pic^0(C)$, then $S$ is representable over $k$.
In particular, if the number of elements $\#(T)$ of $T$ is square-free, 
then $S$ is representable over $k$.
\end{lem}

\begin{proof}
Let $G$ be the free Abelian group generated by
$L,E_1,\ldots,E_r$. We can regard $T$ as a subgroup of 
$G/\langle S\rangle$. If there are infinitely many $l$
such that $T\times {\bf Z}/l{\bf Z}$ is isomorphic to a subgroup of 
$\Pic^0(C)$, then since $\C N_r$ is finite,
we can pick such an $l$ that also has the property that
there is a surjective homomorphism
$\phi :K^\perp \to T\times{\bf Z}/l{\bf Z}$ factoring through the canonical quotient
$K^\perp \to K^\perp/\langle S\rangle$ such that the only elements of
$\C N_r$ in $\ker(\phi)$ are the elements of $\C N_r\cap \langle S\rangle$.
Identifying $T\times {\bf Z}/l{\bf Z}$ with a subgroup of 
$\Pic^0(C)$, we may regard $\phi$ as giving a homomorphism
to $\Pic^0(C)$. If $K^\perp/\langle S\rangle=T$, we proceed as before
using $l=1$. Since $K^\perp$ and $E_1$ generate
$G$ freely, we can extend $\phi$ to a homomorphism $\Phi:G\to \Pic(C)$
by mapping $E_1$ to an arbitrary smooth point $p_1$ of $C$,
and taking $\Phi|_{K^\perp}=\phi$. The images $\Phi(E_i)$ 
now give points $p_i$ of $C$ since $E_i-E_1\in K_X^\perp$, so 
$\Phi(E_i)=\Phi(E_1+(E_i-E_1)) = p_1+\phi(E_i-E_1)$ 
is linearly equivalent to a unique point of $C$, which 
we define to be $p_i$. Blowing up the points $p_i$
gives a morphism $X\to \pr2$ of surfaces 
such that $X$ has no prime divisor
$D$ of self-intersection less than $-2$. Indeed, 
suppose there were such a $D$.
Note that $(-K_X)^2=9-r>0$ and that $-K_X$ is the class of the proper transform of $C$.
Let us denote this proper transform by $C'$. Since 
$-K_X$ is the class of $C'$, a reduced irreducible curve of positive self-intersection, 
we see that $-K_X$ is nef, but by adjunction $D^2<-2$ implies
$-K_X\cdot D<0$, so $D$ cannot be effective.

We want to show that $S=\nneg(X)$. To do this, we will use two facts about $S$.
The first is that if $v\in \langle S\rangle$ has $v^2=-2$
and $v\cdot L\ge 0$, then
$v$ is a nonnegative integer linear combination of elements of $S$.
First we justify this fact. Since $r\le 8$, $K_X^\perp$ is negative definite and even
(i.e., $0\ne v\in K_X^\perp$ implies $v^2$ is negative and even). Suppose 
$v\in \langle S\rangle$ has $v^2=-2$. Write $v=\sum_iP_i-\sum_jN_j$
for elements $P_i$ and $N_j$ of $S$, where $P_i\ne N_j$ for
all $i$ and $j$. If $\sum_iP_i=0$ or $\sum_jN_j=0$, then of course
either $v$ or $-v$ is a nonnegative integer linear combination of elements of $S$.
But $S\subseteq \C N_r\cap K_X^\perp$, and every element of $\C N_r\cap K_X^\perp$
meets $L$ positively, so $-v$ cannot be a nonnegative linear combination of
elements of $S$. I.e., $v$ is a nonnegative integer linear combination of elements of $S$,
as claimed. To finish, we have one case remaining to consider:
suppose $\sum_iP_i\ne0$ and $\sum_jN_j\ne0$.
Then, since $S$ is contained in $K_X^\perp$
and is pairwise nonnegative, we have
$-2=v^2=(\sum_iP_i)^2+(\sum_jN_j)^2-2(\sum_iP_i)(\sum_jN_j)\le -2-2-0$,
which is impossible. 

The second fact about $S$ is that elements of $S$ are effective.
Recall that $-K_X$ is the class of $C'$, which above we noted is nef
of positive self-intersection.
So if $D\in S$, then $(D-C')\cdot C' <0$, hence $h^0(X, D-C')=0$.
Moreover, $h^2(X, D-C')=h^0(X,-D)=0$
since $D\cdot L>0$, and now $h^1(X, D-C')=0$ follows by Riemann-Roch.
It follows that
the restriction morphism $\C O_X(D)\to \C O_{C'}(D)$ is surjective on
global sections, but $h^0(C', \C O_{C'}(D))=1$
since $D\in S\subset \ker(f)$; i.e., $h^0(X, D)=1$
so $D$ is effective.

We now show that $\nneg(X)\subseteq S$.
If $D$ is a prime divisor with $D^2=-2$,
then $-K_X\cdot D=0$ and $D\in \ker(\phi)$, so by construction
$D\in \langle S\rangle$. As shown above, we thus have
$D=\sum_iP_i$ for elements $P_i$ of $S$. 
Since each $P_i$ is effective and $D$ is prime of negative self-intersection, 
$D$ must be one of the $P_i$, hence $D\in S$ so $\nneg(X)\subseteq S$.

Finally, we show that $S\subseteq \nneg(X)$. Say $D\in S$, so $D$ is 
effective. Since $D$ has negative self-intersection, $D$ must meet an element
of $\nneg(X)\subseteq S$ negatively. Since $S$ is pairwise
nonnegative, this element can only be $D$ itself, so 
$D\in \nneg(X)$, hence $S\subseteq\nneg(X)$.
Thus $S=\nneg(X)$ so $S$ is representable.

If $\#(T)$ is square-free, then for any square-free $l$ relatively prime to $\#(T)$,
the group $T\times {\bf Z}/l{\bf Z}$ is cyclic and its order
is square-free. Since any smooth non-supersingular plane cubic $C$
has cyclic torsion subgroups of all square-free orders,
$T\times {\bf Z}/l{\bf Z}$ is isomorphic to a subgroup of
$\Pic^0(C)$.
\end{proof}

We can now prove Theorem \ref{7ptrepthm}:

\begin{proof}
The proof is either by Lemma \ref{Torsionlem} or 
by reduction to previous cases using Lemma \ref{7ptreplem},
so we begin by verifying that every formal configuration type for $r\le 6$ is
representable. This is obvious for $r=1$, and now Lemma \ref{7ptreplem}
applies up through $r=5$. All formal types for $r=6$ also follow 
from $r=5$ by Lemma \ref{7ptreplem} except for type 10, since
$\#(S(6)-S(5))=2$ for type 10 if we take the points in the order
given in our list of types for $r=6$. In this case $\#(T)=2$, however, so
representability follows by Lemma \ref{Torsionlem}. 

Now consider $r=7$. Checking the table of types for $r=7$ we see,
taking the points in the order given in the table, that
Lemma \ref{7ptreplem} applies for all the types except 
21, 23, 24, 28 and 29. For types 21, 28 and 29,
$T$ has order 0, so these types are representable by Lemma \ref{Torsionlem}.
For types 23 and 24, note that up to the general linear group,
we may choose the points $p_1,\ldots,p_7$ to be, respectively,
$(1,0,0)$, $(0,1,0)$, $(0,0,1)$, $(1,1,0)$, $(1,0,1)$,
$(0,1,1)$, and $(1,1,1)$. If we blow these points up to get
$X$, and if $S$ has type 23, then clearly
$S\subset \nneg(X)$. Moreover, the only element of
$\C N_7$ which meets every element of $S$ nonnegatively
is $L-E_4-E_5-E_6$ (i.e., the proper transform of the line through the points
$(1,1,0)$, $(1,0,1)$ and $(0,1,1)$); adding this class to $S$ gives 
type 24. Thus type 23 is representable
if and only if these three points are not collinear,
and type 24 is representable if and only if they are collinear.
But they are collinear if and only if $k$ has characteristic 2.
\end{proof}

For types with $r=8$ points we have:

\begin{thm}\label{8ptrepthm} Let $k$ be an algebraically closed field.
Consider formal configuration types for $r=8$ points.
Types 23, 31, 44, 90, 112, 128, 131 are representable if and only $\Char(k)\neq2$.
Types 46 and 130 are representable if and only if $\Char(k)=2$.
Types 30, 45 and 96 are never representable.
The rest are always representable.
\end{thm}

We will need a version of Lemma \ref{7ptreplem} for $r=8$.

\begin{lem}\label{8ptreplem} Let $k$ be an algebraically closed field,
and let $S$ be a formal configuration type for $r=8$
points. If $S(7)$ is $k$-representable, $\#(S-S(7))\le1$,
and if for each class $D\in \C N_8$ with $D\cdot L=3$ and $D^2 < -1$
there is a $C\in S$ such that $C\ne D$ but $C\cdot D<0$, 
then $S$ is $k$-representable. 
\end{lem}

\begin{proof}
The proof is the same as for Lemma \ref{7ptreplem}.
That proof assumes that $D\in\C N_r$ with $D^2 < -1$ implies
that $D\cdot E_i$ is always either 0 or 1. For $r=8$, this can fail
since, for example, $3L-2E_1-E_2-\cdots-E_8\in \C N_r$.
Thus it is possible a priori that such a class is in
$\nneg(X)$ for the surface obtained in the proof of
Lemma \ref{7ptreplem}. The hypothesis that for each such $D$
there is a $C\in S$ such that $C\cdot D<0$
guarantees that this does not happen.
\end{proof}

We now prove Theorem \ref{8ptrepthm}:

\begin{proof}
For the first 32 types, we have $S\subset K^\perp$.
If $S$ is representable, then $S=\nneg(X)$ for some $X$,
and for this $X$ the class $-K_X$ is nef.
The torsion subgroups $T$, when nonzero, are as follows:
${\bf Z}/2{\bf Z}$ for types  13, 16, 19, 24, 25, and 29, 
$({\bf Z}/2{\bf Z})^2$ for types  23 and 31, 
$({\bf Z}/2{\bf Z})^3$ for type  30, 
${\bf Z}/3{\bf Z}$ for types  27 and 28, 
and 
$({\bf Z}/3{\bf Z})^2$ for type  32. 
Thus, by Lemma \ref{Torsionlem}, possibly except for
types 23, 31, 30 and 32, they are always representable.

For type 32 it turns out that $K^\perp/\langle S\rangle=T$.
This $T$ embeds in $\Pic^0(C)$ for any smooth cubic $C$
if $\Char(k)\neq3$ and for any cuspidal cubic
if $\Char(k)=3$. Thus type 32 is always representable, 
by Lemma \ref{Torsionlem}. 

By Lemma \ref{Torsionlem}, types 23 and 31 are representable if 
$\Char(k)\neq2$. If either of these types $S$
were representable in characteristic 2, then
$S(7)$ would be representable, but $S(7)$ in
either case is the $r=7$ point type
numbered 23 in our table, which is not representable
in characteristic 2. 

For type 30, $S$ is never representable. Suppose it is representable.
Then $-K_X$ is nef, hence $|-K_X|$ is a pencil with no fixed components, by Proposition \ref{8pts}. 
Thus $|-K_X|$ contains an integral divisor; i.e., $X$ is obtained by blowing
up smooth points on a reduced irreducible cubic $C$, where the proper
transform of $C$ is the integral element of $|-K_X|$.
But the points $p_i\in C$ blown up to give $X$ are distinct,
hence the images $E_1-E_i$ for $i\ge1$ under
$\Phi:\Pic(X)\to \Pic(C)$ are distinct. The images of $E_1-E_i$ for $i\le8$
factor through $T$. Since the points $p_i$ are distinct, the image of $T$ in
$\Pic^0(C)$ has order at least 7, whereas $\#(T)$ = 8, so
in fact $\Phi$ gives an injection of $T$ into $\Pic^0(C)$.
But the only cubic curve
whose Picard group contains $({\bf Z}/2{\bf Z})^3$
is the cuspidal cubic, in characteristic 2. 
Thus $\Phi(K^\perp)$ has pure 2-torsion.
We can write $D=3L-E_1-\cdots -E_7-2E_8$
in terms of $K^\perp$ and the elements of $S$ as given in the table
of 8 point types, i.e.,
$D\in \langle S\rangle + 2K^\perp$, hence
$\Phi(D)=0$. Thus $D$ is effective. Since
$D\not\in \langle S\rangle$, this means that
$\nneg(X)\neq S$.

Types 33 through 50 are exactly those for which
$S$ contains a cubic; in particular,
$C'=3L-E_1-\cdots-E_7-2E_8\in S$.
If $C'$ and $D$ are in $\nneg(X)$ and $C'\neq D$, then
$D\cdot E_8=0$, since a check of all elements 
$D\in\C N_8$ shows that if $D\cdot E_8>0$, then
$D\cdot C'<0$. Thus for any type $S$ containing $C'$
we have $S(7)=S-\{C'\}$, and $S$ is representable
if and only if the surface $X$ which represents
$S$ comes from blowing up points $p_1,\ldots, p_8$,
where $p_1,\ldots, p_7$ are smooth points on a 
reduced irreducible singular cubic $C$ with $p_8$
being the singular point of $C$.
Thus, if $S(7)$ is not representable, neither is $S$.
And if $S(7)$ is representable by smooth points on 
a reduced irreducible singular cubic $C$, then
so is $S$. 

Let $T$ be the torsion group for $S(7)$.
For types $S$ numbered 33 through 50 it turns out that
$T$ is 0 except as follows:
it is ${\bf Z}/2{\bf Z}$ for types 41, 44;
$({\bf Z}/2{\bf Z})^2$ for type 45; and
$({\bf Z}/2{\bf Z})^3$ for type 46.
When $T$ is zero, $S(7)$ is (by applying Lemma \ref{Torsionlem})
representable by blowing up smooth points of a cubic with a node,
since in that case $\Pic^0(C)$ is the multiplicative group $k^*$.

As for type 46, by taking 7 smooth points on a cuspidal
cubic in characteristic 2 we can, by Lemma \ref{Torsionlem}, 
represent $S(7)$, hence as observed above, $S$ is representable 
in characteristic 2 for type 46. Since $S(7)$ is not representable
in characteristic not 2, neither is $S$.

Types 41 and 44 are representable by Lemma \ref{Torsionlem} if
the characteristic is not 2. Type 41 is representable also in characteristic
2, since taking $C$ to be a cuspidal cubic, the
homomorphism $K^\perp\to\Pic^0(C)$ factors through
$K^\perp/(\langle S(7)\rangle +2K^\perp)=({\bf Z}/2{\bf Z})^4$,
but by explicitly checking the map,
no element of $\C N_7$ is in the kernel except those of $S(7)$.
On the other hand, type 44 is not representable in characteristic 2, since
if it were we would need to blow up 7 smooth points on a cuspidal cubic $C$,
but then $L_{123},L_{145},L_{257}\in S(7)$ and
$L_{347}$ and $L_{356}$ would be in the kernel of $K^\perp\to \Pic^0(C)$ even though
neither is in $S(7)$; e.g.,
$L_{347}\equiv L_{123}+L_{145}+L_{257}\hbox{ (mod 2)}$,
where $L_{i_1i_2i_3}=L-E_{i_1}-E_{i_2}-E_{i_3}$.

Type 45 is not representable in characteristic 2,
since $S(7)$ is not. It is not representable in any other characteristic
either. If it were, we would need to find smooth points on a singular
cubic $C$ such that for the induced homomorphism $\Pic(X_7)\to \Pic(C)$,
where $X_7$ is the blow-up of the first 7 points and $p_8$ is the singular
point of $C$, the only $(-2)$-classes in the kernel are the elements
of $S(7)$. But $K_{X_7}^\perp/\langle S(7)\rangle$ is 
$({\bf Z}/2{\bf Z})^2\oplus {\bf Z}$. Since the characteristic is not 2
but $C$ is singular,
the torsion of the image of this quotient in $\Pic(C)$ is cyclic. Thus some element of the torsion
subgroup of $K_{X_7}^\perp/\langle S(7)\rangle$ must map to 0. 
But there are three such elements, and if $x$ is any one of them, an explicit
check shows that $\langle S(7), x\rangle$ contains $(-2)$-classes not in
$S(7)$, hence $X_7$ would have $\nneg(X_7)\ne S(7)$.

For types 51 through 96, $S\subset K^\perp$, so we can apply
Lemma \ref{Torsionlem}. The torsion subgroup $T$ of $K^\perp/\langle S\rangle$
is zero except as follows:
${\bf Z}/2{\bf Z}$ for types 58, 60, 78, 82, 86, 89, 94;
$({\bf Z}/2{\bf Z})^2$ for types 90, 95;
$({\bf Z}/2{\bf Z})^3$ for type 96; and
${\bf Z}/3{\bf Z}$ for type 62.
Thus all are always representable except possibly
types 90, 95 and 96. Types 90 and 95 are representable
(taking points on a smooth non-supersingular cubic, by Lemma \ref{Torsionlem})
except possibly in characteristic 2. 
Type 90 is not representable in characteristic 2.
We cannot take smooth points on a cuspidal cubic $C$, since
then the kernel of $\Pic(X)\to \Pic(C)$ contains
$\langle S\rangle+2K^\perp$, which by explicit check
contains $(-2)$-classes other than those in $S$.
If $S$ is representable by choosing smooth points
on either a nodal cubic or a smooth cubic, then the 2-torsion
is at most ${\bf Z}/2{\bf Z}$, so the map
$({\bf Z}/2{\bf Z})^2=T\to \Pic^0(C)$
induced by $\Pic(X)\to \Pic(C)$
must kill some of the 2-torsion of $T$. But by brute force
check, every 2-torsion element of $T$ is the image of a $(-2)$-class,
hence killing any 2-torsion makes $\nneg(X)\neq S$.
Type 95 is representable in characteristic 2.
We cannot take smooth points on a cuspidal cubic $C$, since
then the kernel of $\Pic(X)\to \Pic(C)$ contains
$\langle S\rangle+2K^\perp$, which by explicit check
contains $(-2)$-classes other than those in $S$, but
there is an element $x$ of $K^\perp$ such that
$x$ maps to a 2-torsion element of $K^\perp/\langle S\rangle$
but such that the only $(-2)$-classes in 
$\langle S, x\rangle$ are those in $S$.
Thus $K^\perp/\langle S, x\rangle$ has torsion group $T'$
which embeds in $\Pic(C)$ for either a nodal or smooth
and non-supersingular $C$ in characteristic 2.
Thus, as in Lemma \ref{Torsionlem}, $S$ is representable
in characteristic 2.

Type 96 is never representable. If the characteristic is not 2,
we must kill some of the 2-torsion, but by the same method
as in the case of type 90, doing so introduces
extra $(-2)$-classes. Thus the only hope for representability
is in characteristic 2, with points on a cuspidal cubic.
But $\langle S\rangle+2K^\perp$ by explicit check
contains $(-2)$-classes other than those in $S$,
so even this does not work.

We need to treat types
111, 112, 119, 121, 126, 128, 129, 130, 131
specially. For type 111, a brute force check shows
that there is no $(-2)$-class
which meets the elements of $S$ nonnegatively.
Thus if we can choose points such that
$S\subset \nneg(X)$, then $S=\nneg(X)$. But in this case
it is easy to check that by choosing our points
from among the intersections of three general lines and 
a general conic we do indeed get $S\subset \nneg(X)$.
For type 112, perform a quadratic transformation centered at
the points c, d, and g. This transforms type 112 into type 128.
Thus the one is representable if and only if the other is.
(Alternatively, we can regard types 112 and 128 as giving the 
same surface $X$, but with respect to different morphisms $X\to \pr2$.)
Since for type 128, $S(7)$ is not representable in characteristic 2,
neither is type 128 (nor 112). But by Lemma \ref{8ptreplem},
type 128 (and hence 112) is representable if the characteristic is not 2.
Type 119 can be handled in the same way that 111 was.
Type 121 can be handled by applying a quadratic transformation
centered at the points a, b and c. This transforms 121 into type 144,
which can be handled by Lemma \ref{8ptreplem}. 
Type 126 can be handled as was 111, choosing the 8 points from among
the 10 intersections of five general lines (thus it is easy to ensure that
$S\subset \nneg(X)$, but a brute force check shows
that there is no $(-2)$-class
which meets the elements of $S$ nonnegatively.)
For type 129, since all of the elements of $S$ come from
lines through either the point a or the point b, and since e is on the line through
a and d, as long as the points are distinct, e cannot be on any
of the lines  but the line through a and d. Likewise, h can only be on the line
containing c, d and g and the line containing b and e.
Thus if we choose the points so that e is general but so that
$S(7)$ is representable (which we can, since $S(7)$ is the 7
point type numbered 22), then $S\subset \nneg(X)$. 
The only possible additional element of $\nneg(X)$ besides
$S$ allowed by pairwise nonnegativity would come from a line
through c, e and f. But since e is general, e cannot be on the line
through c and f. Thus $S$ is representable.
For type 130, $S(7)$ is representable if and only if the characteristic
is 2. Thus 130 is not representable if the characteristic is not 2,
and by Lemma \ref{8ptreplem}, it is representable if the characteristic
is 2. For type 131, we may choose coordinates such that
a is the point $(0,0,1)$, c is $(0,1,0)$, g is $(1,0,0)$,
e is $(1,1,1)$, hence $b$ is $(0,1,1)$, and $f$ is $(1,0,1)$.
Now d is forced to be the point $(1,1,-2)$. Since in characteristic 2
this is a point on the line through c and g, but $S$ does not allow the points
c, d and g to be collinear, we see $S$ is not representable in characteristic 2.
But if the  characteristic is not 2, then we can check explicitly that
$S\subset \nneg(X)$. No other additional element of $\nneg(X)$ besides
$S$ is allowed by pairwise nonnegativity,
so $S=\nneg(X)$ here.

The remaining types are handled by Lemma \ref{8ptreplem}
in a way similar to what was done applying
Lemma \ref{7ptreplem} in the proof of Theorem  \ref{7ptrepthm}.
For example, for type 97, the type $S'$ given
by the points a, b, c, d, e, g, h is the 7 point type 29.
It is representable, and $\#(S-S')=0$, so 
type 97 is representable, by  Lemma \ref{8ptreplem}.
\end{proof}

\begin{example}\label{infmanytypes} \rm
Here we show briefly that there are infinitely many configuration types among sets of 
$r$ points, for each $r\ge9$. In fact, it is clear by the definition of configuration type
that if there are infinitely many types for $r=9$, then there are infinitely many for all $r>9$.
So pick points $p_1,\ldots,p_9$ on a smooth cubic curve $C'$. 
Let $Z=p_1+\cdots+p_9$, and
let $X$ be the surface obtained by blowing up the points,
and let $C$ be the proper transform of $C'$ on $X$; note that the class of $C$ is $-K_X$.
Using the group law on the cubic, it is not hard to see that among all choices of the points $p_i$
there arise infinitely many different positive integers $m_Z$ such that $m_Z$ is
the least positive integer for which the restriction of $m_ZK_X$ to $C$ is trivial (as a line bundle).
By results of \cite{ARS}, it follows that $h_{m_ZZ}(3i)$ is $\binom{3i+2}{2}-1$ for $i<m_Z$,
while $h_{m_ZZ}(3m_Z)=\binom{3i+2}{2}-2$. I.e., there are infinitely many
configuration types of 9 points.
\end{example}

\begin{rem}\label{sympowersrem} \rm 
It is possible for $Z$ and $Z'$ to have different configuration types
but nonetheless for $mZ$ and $mZ'$ to have the same Hilbert functions for all $m$.
In this situation it is convenient to say that $Z$ and $Z'$ have the same
{\it uniform} configuration type. For example, the 6 point types 8 and 11
have the same uniform type, since they are both a complete intersection of a conic and a cubic.
We also note that in order for two nonequivalent
types to have the same uniform type, it need not be true that they be complete intersections.
In particular, adding a general seventh point to the 6 point types 8 and 11
gives the 7 point types 11 and 26, neither of which is a complete intersection
but for each $m\ge 1$ the Hilbert function of $mZ$ is the same whether $Z$ has type
11 or 26. Type 26 consists of 6 points on an irreducible conic with a general seventh point.
Type 11 consists of two sets of three collinear points with a seventh general point. 
For type 11 it turns out that 
whenever one of the lines through 3 collinear points is a fixed component for forms of degree
$t$ vanishing on $mZ$, the other line through 3 collinear points is, by symmetry,
also a fixed component. Thus the two lines are always taken together, so things 
for $mZ$ when $Z$ has type 11 behave the same as when $Z$ has type 26,
where the two lines are replaced by an irreducible conic.
Finally, it is interesting to mention that one can find reduced finite subschemes $Z$ and $Z'$
of the plane such that $2Z$ and $2Z'$ have the same Hilbert function, but where
the Hilbert functions of $Z$ and $Z'$ are different; see Example 7.1 of \cite{GMS}.
\end{rem}

\begin{rem}\label{oraclerem} \rm 
Because there are only finitely many configuration types of $r$ points
for each $r\le 8$, it follows that there is a number $N_r$ such that
if one knows the Hilbert function of $mZ$ (i.e., of $(I(Z))^{(m)}$) for each $m\le N_r$
for some reduced scheme $Z$ consisting of $r$ points in the plane,
then one can deduce the uniform configuration type of $Z$ and hence
the Hilbert function of $mZ$ for all $m>0$. For example, by examining
the Hilbert function of $p_1+p_2+p_3$ one can tell if the points are collinear or not,
and hence $N_3=1$. By checking Hilbert functions for each type, and by straightforward
arguments to show which different types have the same uniform type, 
we determined that $N_r=1$ for $r\le3$, $N_4=2$, 
$N_r=3$ for $r=5,6$ and $N_7=7$. We have not bothered to
determine exactly which 8 point types have the same uniform types, 
and so we do not know for sure
what the value of $N_8$ is, but it is not less than 16,
and we suspect that it is exactly 16. As an interesting
sidelight, it turns out in fact for each $r\le7$ that the Hilbert function
of $N_rZ$ alone already determines the Hilbert function
of $mZ$ for all $m$, since for $r\le7$, the Hilbert function
of $N_rZ$ distinguishes the uniform type. Thus if $I$ is the ideal of a reduced set of $r\le7$
points, then the Hilbert function of $I^{(N_r)}$ determines the Hilbert functions
of $I^{(m)}$ for all $m>0$. (For the case of $r=8$, the least $N$
for which the Hilbert function of $I^{(N)}$ could by itself determine
the Hilbert functions of $I^{(m)}$ for all $m>0$ is $N=22$,
and we suspect that $N=22$ in fact works.)
\end{rem}
 
\vfil\eject
 
\section{The Tables}\label{ourtables}

\begin{table}[ht]
{\small 
\[
\begin{array}{rlcrlcrl}
\hbox{No.} & \hbox{Type} & \hbox to.3in{\hfil} & \hbox{No.} & \hbox{Type} &  \hbox to.3in{\hfil}&\hbox{No.} & \hbox{Type} \\
\hline
1. & \hbox to.35in{\hfil $\emptyset$\hfil} & & 5. & \hbox{1: abc, ade}& & 9. & \hbox{1: abc, ade, bdf} \\
2. & \hbox{1: abc} & & 6. & \hbox{1: abcdef} & & 10. & \hbox{1: abc, ade, bdf, cef} \\
3. & \hbox{1: abcd} & & 7. & \hbox{1: abcd, aef} & & 11. & \hbox{2: abcdef} \\
4. & \hbox{1: abcde} & & 8. & \hbox{1: abc, def} & & & \\
\end{array}
\]
}
\caption[\hskip.1in Configuration types for $r\leq6$ points]{Configuration types for $r\leq6$ points}\label{6pttableA}
\end{table}

\begin{table}[ht]
\ \vskip-.25in
{\small 
\[
\begin{array}{ccllll}
r & m & \hbox{Type(s)} & h_Z & F_0 & F_1 \\
\hline
1 & 1 & 1 & 1 & 1^2 & 2^1 \\
1 & 2 & 1 & 1, 3 & 2^3 & 3^2 \\
\hline
2 & 1 & 1 & 1, 2 & 1^1, 2^1  & 3^1 \\
2 & 2 & 1 & 1, 3, 5, 6 & 2^1, 3^1, 4^1 & 4^1, 5^1 \\
\hline
3 & 1 & 1 & 1, 3 & 2^3 & 3^2 \\
3 & 2 & 1 & 1, 3, 6, 9 & 3^1, 4^3 & 5^3 \\
\hline
3 & 1 & 2 & 1, 2, 3 & 1^1, 3^1 & 4^1 \\
3 & 2 & 2 & 1, 3, 5, 7, 8, 9 & 2^1, 4^1, 6^1 & 5^1, 7^1 \\
\hline
4 & 1 & 1 & 1, 3, 4 & 2^2 & 4^1 \\
4 & 2 & 1 & 1, 3, 6, 10, 12 & 4^3 & 6^2 \\
\hline
4 & 1 & 2 & 1, 3, 4 & 2^2, 3^1 & 3^1, 4^1 \\
4 & 2 & 2 & 1, 3, 6, 10, 11, 12 & 4^4, 6^1 & 5^3, 7^1 \\
\hline
4 & 1 & 3 & 1, 2, 3, 4 & 1^1, 4^1 & 1^5 \\
4 & 2 & 3 & 1, 3, 5, 7, 9, 10, 11, 12 & 2^1, 5^1, 8^1 & 6^1, 9^1 \\
\hline
5 & 1 & 1, 2 & 1, 3, 5 & 2^1, 3^2 & 4^2 \\
5 & 2 & 1 & 1, 3, 6, 10, 14, 15 & 4^1, 5^3 & 6^2, 7^1 \\
5 & 2 & 2 & 1, 3, 6, 10, 14, 15 & 4^1, 5^3, 6^1 & 6^3, 7^1 \\
\hline
5 & 1 & 3 & 1, 3, 4, 5 & 2^2, 4^1 & 3^1, 5^1 \\
5 & 2 & 3 & 1, 3, 6, 10, 12, 13, 14, 15 & 4^3, 5^1, 8^1 & 5^2, 6^1, 9^1 \\
\hline
5 & 1 & 4 & 1, 2, 3, 4, 5 & 1^1, 5^1 & 6^1 \\
5 & 2 & 4 & 1, 3, 5, 7, 9, 11, 12, 13, 14, 15 & 2^1, 6^1, 10^1 & 7^1, 11^1 \\
\hline
5 & 1 & 5 & 1, 3, 5 & 2^1, 3^2 & 4^2 \\
5 & 2 & 5 & 1, 3, 6, 10, 13, 15 & 4^2, 6^2 & 6^1, 7^2 \\
\hline
6 & 1 & 1, 2, 5, 9, 10 & 1, 3, 6 & 3^4 & 4^3 \\
6 & 2 & 1, 2 & 1, 3, 6, 10, 15, 18 & 5^3, 6^1 & 7^3 \\
6 & 2 & 5 & 1, 3, 6, 10, 15, 18 & 5^3, 6^2 & 6^1, 7^3 \\
6 & 2 & 9 & 1, 3, 6, 10, 15, 18 & 5^3, 6^3 & 6^2, 7^3 \\
6 & 2 & 10 & 1, 3, 6, 10, 14, 18 & 4^1, 6^4 & 7^4 \\
\hline
6 & 1 & 3, 7 & 1, 3, 5, 6 & 2^1, 3^1, 4^1 & 4^1, 5^1 \\
6 & 1 & 8, 11 & 1, 3, 5, 6 & 2^1, 3^1 & 5^1 \\
6 & 2 & 3 & 1, 3, 6, 10, 14, 16, 17, 18 & 4^1, 5^2, 8^1 & 6^1, 7^1, 9^1 \\
6 & 2 & 7 & 1, 3, 6, 10, 14, 16, 17, 18 & 4^1, 5^2, 6^1, 8^1 & 6^2, 7^1, 9^1 \\
6 & 2 & 8, 11 & 1, 3, 6, 10, 14, 17, 18 & 4^1, 5^1, 6^1 & 7^1, 8^1 \\
\hline
6 & 1 & 4 & 1, 3, 4, 5, 6 & 2^2, 5^1 & 3^1, 6^1 \\
6 & 2 & 4 & 1, 3, 6, 10, 12, 14, 15, 16, 17, 18 & 4^3, 6^1, 10^1 & 5^2, 7^1, 11^1 \\
\hline
6 & 1 & 6 & 1, 2, 3, 4, 5, 6 & 1^1, 6^1 & 7^1 \\
6 & 2 & 6 & 1, 3, 5, 7, 9, 11, 13, 14, 15, 16, 17, 18 & 2^1, 7^1, 12^1 & 8^1, 13^1 \\
\end{array}
\]
}
\caption[\hskip.1in Hilbert functions by configuration type for $r\le6$ points]{Hilbert functions by configuration type for $r\le6$ points}\label{6pttableB}
\end{table}

\begin{table}[ht]
{\small 
\[
\begin{array}{clccl}
\hbox{No.} & \hbox{Type} & \hbox to .5in{\hfil} & \hbox{No.} & \hbox{Type}\\
\hline
1 & \hbox{empty} & & 16 & \hbox{1: abc, ade, cef} \\
2 & \hbox{1: abcdefg} & & 17 & \hbox{1: abcg, ade, bdf, cef} \\
3 & \hbox{1: abcdef} & & 18 & \hbox{1: abc, ade, bdf, ceg} \\
4 & \hbox{1: abcde} & & 19 & \hbox{1: abc, ade, cef, afg} \\
5 & \hbox{1: abcd} & & 20 & \hbox{1: abc, adf, cef, bde} \\
6 & \hbox{1: abc} & & 21 & \hbox{1: abc, def, adg, beg, cfg} \\
7 & \hbox{1: abcde, afg} & & 22 & \hbox{1: abc, adf, cef, bde, aeg} \\
8 & \hbox{1: abcd, efg} & & 23 & \hbox{1: abc, adf, cef, bde, aeg, cdg} \\
9 & \hbox{1: abcd, defg} & & 24 & \hbox{1: abc, adf, cef, bde, aeg, cdg, bfg} \\
10 & \hbox{1: abcd, def} & & 25 & \hbox{2: abcdefg} \\
11 & \hbox{1: abc, def} & & 26 & \hbox{2: abcdef} \\
12 & \hbox{1: abc, ade} & & 27 & \hbox{1: abg; 2: abcdef} \\
13 & \hbox{1: abcd, def, ceg} & & 28 & \hbox{1: abg, cdg; 2: abcdef} \\
14 & \hbox{1: abc, def, adg} & & 29 & \hbox{1: abg, cdg, efg; 2: abcdef} \\
15 & \hbox{1: abc, ade, afg} & & & 
\end{array}
\]
}
\caption[\hskip.1in Seven point configuration types]{Seven point configuration types}\label{7pttableA}
\end{table}

\begin{table}[ht]
{\small 
\[
\begin{array}{ccllll}
r & m & \hbox{Type(s)} & h_Z \\
\hline
7 & 1 & \hbox{8, 9, 25} & \hbox{1, 3, 5, 7} \\
7 & 2 & \hbox{8, 25} & \hbox{1, 3, 6, 10, 14, 18, 20, 21} \\
7 & 2 & 9 & \hbox{1, 3, 6, 10, 14, 17, 19, 21} \\
\hline
7 & 1 & \hbox{1, 5, 6, 10, 11, 12, 13, $\ldots$, 24, 26, $\ldots$, 29} & \hbox{1, 3, 6, 7} \\
7 & 2 & \hbox{11, 14, 18, 21, 26, 27, 28, 29} & \hbox{1, 3, 6, 10, 15, 20, 21} \\
7 & 2 & \hbox{5, 10, 13, 17} & \hbox{1, 3, 6, 10, 15, 19, 20, 21} \\
7 & 2 & \hbox{1, 6, 12, 15, 16, 19, 20, 22, 23, 24} & \hbox{1, 3, 6, 10, 15, 21} \\
\hline
7 & 1 & 2 & \hbox{1, 2, 3, 4, 5, 6, 7} \\
7 & 2 & 2 & \hbox{1, 3, 5, 7, 9, 11, 13, 15, 16, 17, 18, 19, 20, 21} \\
\hline
7 & 1 & 3 & \hbox{1, 3, 4, 5, 6, 7} \\
7 & 2 & 3 & \hbox{1, 3, 6, 10, 12, 14, 16, 17, 18, 19, 20, 21} \\
\hline
7 & 1 & \hbox{4, 7} & \hbox{1, 3, 5, 6, 7} \\
7 & 2 & \hbox{4, 7} & \hbox{1, 3, 6, 10, 14, 17, 18, 19, 20, 21} \\
\end{array}
\]
}
\caption[\hskip.1in Hilbert functions by configuration type for $r=7$ points]{Hilbert functions by configuration type for $r=7$ points}\label{7pttableB}
\end{table}

\vbox to2in{\vfil}

\newpage
{\tiny 
\[
\begin{array}{clcl}
\hbox{No.} & \hbox{Type} &  \hbox{No.} & \hbox{Type}\\
\hline
1 & \hbox{empty} & 51 & \hbox{2: abcdef} \\
2 & \hbox{1: abc} & 52 & \hbox{2: abcdef, abcdgh} \\
3 & \hbox{1: abc, def} & 53 & \hbox{2: abcdef, abcdgh, abefgh} \\
4 & \hbox{1: abc, ade} & 54 & \hbox{2: abcdef, abcdgh, abefgh, cdefgh} \\
5 & \hbox{1: abc, ade, afg} & 55 & \hbox{1: abc, ade, fgh, 2: bcdegh} \\
6 & \hbox{1: abc, ade, bdf} & 56 & \hbox{1: abc, ade, bdf, cgh, 2: abefgh} \\
7 & \hbox{1: abc, ade, bfg} & 57 & \hbox{1: abc, ade, afg, bdf, 2: cdefgh} \\
8 & \hbox{1: abc, ade, fgh} & 58 & \hbox{1: abc, ade, afg, bdf, beg, 2: cdefgh} \\
9 & \hbox{1: abc, ade, bdf, cgh} & 59 & \hbox{1: abc, ade, afg, bdf, beh, 2: cdefgh} \\
10 & \hbox{1: abc, ade, bdf, ceg} & 60 & \hbox{1: abc, ade, afg, bdh, ceh, 2: bcdefg} \\
11 & \hbox{1: abc, ade, bdf, cef} & 61 & \hbox{1: abc, ade, afg, bdh, cfh, 2: bcdefg} \\
12 & \hbox{1: abc, ade, bfg, dfh} & 62 & \hbox{1: abc, ade, afg, bdh, cfh, egh, 2: bcdefg} \\
13 & \hbox{1: abc, ade, afg, bdf} & 63 & \hbox{1: abc, 2: cdefgh} \\
14 & \hbox{1: abc, ade, afg, bdh} & 64 & \hbox{1: abc, ade, 2: bcdegh} \\
15 & \hbox{1: abc, ade, afg, bdf, ceg} & 65 & \hbox{1: abc, ade, afg, 2: bcdefg} \\
16 & \hbox{1: abc, ade, afg, bdf, beg} & 66 & \hbox{1: abc, ade, bdf, 2: cdefgh} \\
17 & \hbox{1: abc, ade, afg, bdf, ceh} & 67 & \hbox{1: abc, ade, bfg, 2: cdefgh} \\
18 & \hbox{1: abc, ade, afg, bdf, beh} & 68 & \hbox{1: abc, ade, afg, bdh, 2: cdefgh} \\
19 & \hbox{1: abc, ade, afg, bdh, ceh} & 69 & \hbox{1: abc, 2: bcdefg} \\
20 & \hbox{1: abc, ade, afg, bdh, cfh} & 70 & \hbox{1: abc, ade, 2: cdefgh} \\
21 & \hbox{1: abc, ade, bdf, cgh, efg} & 71 & \hbox{1: abc, ade, afg, 2: cdefgh} \\
22 & \hbox{1: abc, ade, afg, bdf, ceg, beh} & 72 & \hbox{1: abc, ade, bdf, 2: bcdegh} \\
23 & \hbox{1: abc, ade, afg, bdf, beg, cdg} & 73 & \hbox{1: abc, ade, bfg, 2: bcdegh} \\
24 & \hbox{1: abc, ade, afg, bdf, beg, dgh} & 74 & \hbox{1: abc, ade, afg, bdh, 2: bcdefg} \\
25 & \hbox{1: abc, ade, afg, bdf, beg, cdh} & 75 & \hbox{1: abc, ade, 2: acefgh} \\
26 & \hbox{1: abc, ade, afg, bdf, ceh, bgh} & 76 & \hbox{1: abc, def, 2: bcefgh} \\
27 & \hbox{1: abc, ade, afg, bdh, cfh, egh} & 77 & \hbox{1: abc, ade, bdf, ceg, 2: acdfgh} \\
28 & \hbox{1: abc, ade, afg, bdf, ceg, beh, cfh} & 78 & \hbox{1: abc, ade, bdf, cef, 2: bcdegh} \\
29 & \hbox{1: abc, ade, afg, bdf, ceg, beh, cdh} & 79 & \hbox{1: abc, ade, bfg, dfh, 2: bcdegh} \\
30 & \hbox{1: abc, ade, afg, bdf, beg, cdg, cef} & 80 & \hbox{1: abc, 2: cdefgh, abefgh} \\
31 & \hbox{1: abc, ade, afg, bdf, beg, cdg, ceh} & 81 & \hbox{1: abc, ade, 2: bcdegh, acefgh} \\
32 & \hbox{1: abc, ade, afg, bdf, ceg, beh, cfh, dgh} & 82 & \hbox{1: abc, ade, bdf, 2: bcdegh, acdfgh} \\
33 & \hbox{3: abcdefgh} & 83 & \hbox{1: abc, 2: acdefg, abefgh} \\
34 & \hbox{1: abc, 3: abcdefgh} & 84 & \hbox{1: abc, ade, 2: bdefgh, acefgh} \\
35 & \hbox{1: abc, def, 3: abcdefgh} & 85 & \hbox{1: abc, ade, bdf, 2: cdefgh, abefgh} \\
36 & \hbox{1: abc, ade, 3: abcdefgh} & 86 & \hbox{1: abc, ade, 2: acefgh, abdfgh} \\
37 & \hbox{1: abc, ade, afg, 3: abcdefgh} & 87 & \hbox{1: abc, def, 2: bcefgh, acdfgh} \\
38 & \hbox{1: abc, ade, bdf, 3: abcdefgh} & 88 & \hbox{1: abc, ade, bfg, 2: bcdegh, acefgh} \\
39 & \hbox{1: abc, ade, bfg, 3: abcdefgh} & 89 & \hbox{1: abc, ade, bdf, ceg, 2: acdfgh, abefgh} \\
40 & \hbox{1: abc, ade, bdf, ceg, 3: abcdefgh} & 90 & \hbox{1: abc, ade, bdf, cef, 2: bcdegh, acdfgh} \\
41 & \hbox{1: abc, ade, bdf, cef, 3: abcdefgh} & 91 & \hbox{1: abc, ade, bfg, dfh, 2: bcdegh, acefgh} \\
42 & \hbox{1: abc, ade, afg, bdf, 3: abcdefgh} & 92 & \hbox{1: abc, 2: bcdfgh, acdefg, abefgh} \\
43 & \hbox{1: abc, ade, afg, bdf, ceg, 3: abcdefgh} & 93 & \hbox{1: abc, def, 2: bcefgh, acdfgh, abdegh} \\
44 & \hbox{1: abc, ade, afg, bdf, beg, 3: abcdefgh} & 94 & \hbox{1: abc, ade, 2: bcdegh, acefgh, abdfgh} \\
45 & \hbox{1: abc, ade, afg, bdf, beg, cdg, 3: abcdefgh} & 95 & \hbox{1: abc, ade, bdf, 2: bcdegh, acdfgh, abefgh} \\
46 & \hbox{1: abc, ade, afg, bdf, beg, cdg, cef, 3: abcdefgh} & 96 & \hbox{1: abc, ade, bdf, cef, 2: bcdegh, acdfgh, abefgh} \\
47 & \hbox{2: abcdef, 3: abcdefgh} & 97 & \hbox{1: abc, ade, afgh, 2: bcdegh} \\
48 & \hbox{1: abc, 2: bcdefg, 3: abcdefgh} & 98 & \hbox{1: abc, defg} \\
49 & \hbox{1: abc, ade, 2: bcdefg, 3: abcdefgh} & 99 & \hbox{1: abc, ade, afgh} \\
50 & \hbox{1: abc, ade, afg, 2: bcdefg, 3: abcdefgh} & 100 & \hbox{1: abc, ade, bdf, afgh} \\
\end{array}
\]
}
\begin{table}[ht]
\caption[\hskip.1in Eight point configuration types (part 1)]{Eight point configuration types (part 1)}\label{8pttableA1}
\vskip-.3in
\end{table}

\newpage

\begin{table}[ht]
{\small 
\[
\begin{array}{clccl}
\hbox{No.} & \hbox{Type} & \hbox to .2in{\hfil} & \hbox{No.} & \hbox{Type}\\
\hline
101 & \hbox{1: abc, ade, afg, bdf, cegh} & & 124 & \hbox{1: abc, ade, bdf, ceg, afgh} \\
102 & \hbox{1: abc, adeh, 2: bcdefg} & & 125 & \hbox{1: abc, ade, bdf, cef, afgh} \\
103 & \hbox{1: abc, ade, bdf, afgh, 2: bcdegh} & & 126 & \hbox{1: abc, ade, bfg, dfh, cegh} \\
104 & \hbox{1: abc, adef} & & 127 & \hbox{1: abc, ade, afg, bdh, cefh} \\
105 & \hbox{1: abc, ade, bdfg} & & 128 & \hbox{1: abc, ade, afg, bdf, beg, cdgh} \\
106 & \hbox{1: abc, ade, bdf, cegh} & & 129 & \hbox{1: abc, ade, afg, bdf, beh, cdgh} \\
107 & \hbox{1: abc, ade, afg, bdf, begh} & & 130 & \hbox{1: abc, ade, afg, bdf, beg, cdg, cefh} \\
108 & \hbox{1: abc, ade, bdfg, 2: acefgh} & & 131 & \hbox{1: abc, ade, afg, bdf, beg, dgh, cefh} \\
109 & \hbox{1: abc, ade, bfgh} & & 132 & \hbox{1: abcde} \\
110 & \hbox{1: abc, ade, bdf, cefg} & & 133 & \hbox{1: abcde, afgh} \\
111 & \hbox{1: abc, def, adgh, 2: bcefgh} & & 134 & \hbox{1: abcde, fgh} \\
112 & \hbox{1: abc, ade, bdf, cef, afgh, 2: bcdegh} & & 135 & \hbox{1: abcde, afg} \\
113 & \hbox{1: abcd} & & 136 & \hbox{1: abcde, afg, bfh} \\
114 & \hbox{1: abc, ade, afg, bdfh} & & 137 & \hbox{1: abcde, afg, bfh, cgh} \\
115 & \hbox{1: abcd, efgh} & & 138 & \hbox{1: abcdef} \\
116 & \hbox{1: abcd, aefg} & & 139 & \hbox{1: abcdef, agh} \\
117 & \hbox{1: abc, adef, bdgh} & & 140 & \hbox{1: abcdefg} \\
118 & \hbox{1: abc, ade, bdfg, cefh} & & 141 & \hbox{1: abcdefgh} \\
119 & \hbox{1: abc, ade, afg, bdfh, cegh} & & 142 & \hbox{2: abcdefg} \\
120 & \hbox{1: abgh, 2: abcdef} & & 143 & \hbox{1: abh, 2: abcdefg} \\
121 & \hbox{1: efgh, 2: abcdef, abcdgh} & & 144 & \hbox{1: abh, cdh, 2: abcdefg} \\
122 & \hbox{1: abc, def, adgh} & & 145 & \hbox{1: abh, cdh, efh, 2: abcdefg} \\
123 & \hbox{1: abc, ade, bfg, cdfh} & & 146 & \hbox{2: abcdefgh}
\end{array}
\]
}
\caption[\hskip.1in Eight point configuration types (part 2)]{Eight point configuration types (part 2)}\label{8pttableA2}
\end{table}

\begin{table}[ht]
{\small 
\[
\begin{array}{ccllll}
r & m & \hbox{Type(s)} & h_Z \\
\hline
8 & 1 & \hbox{1, $\ldots$, 114, 116, $\ldots$, 131, 142, \dots, 145} & \hbox{1, 3, 6, 8} \\
8 & 2 & \hbox{1, $\ldots$, 96} & \hbox{1, 3, 6, 10, 15, 21, 24} \\
8 & 2 & \hbox{97, $\ldots$, 110, 112, 113, 114, 120,} & \\
  &  &  \hbox{122, $\ldots$, 125, 127, $\ldots$, 131, 142, $\ldots$, 145} & \hbox{1, 3, 6, 10, 15, 21, 23, 24} \\
8 & 2 & \hbox{111, 121, 126} & \hbox{1, 3, 6, 10, 15, 20, 23, 24} \\
8 & 2 & \hbox{116, 117, 118, 119} & \hbox{1, 3, 6, 10, 15, 20, 22, 24} \\
\hline
8 & 1 & \hbox{115, 133, 134, 146} & \hbox{1, 3, 5, 7, 8} \\
8 & 2 & \hbox{115, 146} & \hbox{1, 3, 6, 10, 14, 18, 21, 23, 24} \\
8 & 2 & 133 & \hbox{1, 3, 6, 10, 14, 18, 20, 22, 23, 24} \\
8 & 2 & 134 & \hbox{1, 3, 6, 10, 14, 18, 21, 22, 23, 24} \\
\hline
8 & 1 & \hbox{132, 135, 136, 137} & \hbox{1, 3, 6, 7, 8} \\
8 & 2 & \hbox{132, 135, 136, 137} & \hbox{1, 3, 6, 10, 15, 20, 21, 22, 23, 24} \\
\hline
8 & 1 & \hbox{138, 139} & \hbox{1, 3, 5, 6, 7, 8} \\
8 & 2 & \hbox{138, 139} & \hbox{1, 3, 6, 10, 14, 17, 19, 20, 21, 22, 23, 24} \\
\hline
8 & 1 & 140 & \hbox{1, 3, 4, 5, 6, 7, 8} \\
8 & 2 & 140 & \hbox{1, 3, 6, 10, 12, 14, 16, 18, 19, 20, 21, 22, 23, 24} \\
\hline
8 & 1 & 141 & \hbox{1, 2, 3, 4, 5, 6, 7, 8} \\
8 & 2 & 141 & \hbox{1, 3, 5, 7, 9, 11, 13, 15, 17, 18, 19, 20, 21, 22, 23, 24} \\
\end{array}
\]
}
\caption[\hskip.1in Hilbert functions by configuration type for $r=8$ points]{Hilbert functions by configuration type for $r=8$ points}\label{8pttableB}
\end{table}

\end{document}